\documentclass{amsart}
\usepackage{array}
\usepackage{multirow}
\usepackage{graphicx}
\usepackage{epsfig}
\usepackage{color,amsfonts,amsmath}
\usepackage{amssymb}
\usepackage{latexsym}
\usepackage{epsf}
\usepackage{graphicx,color,graphics}
\usepackage{amsmath,amssymb}
\usepackage{dsfont}
\usepackage[latin1]{inputenc}
\newtheorem{theorem}{Theorem}[section]

\newtheorem{lemma}[theorem]{Lemma}

\newtheorem{remark}{Remark}

\numberwithin{equation}{section}

\begin{document}

\title[Cauchy Rao-Nakra sandwich beam]{Some $L^q (\mathbb{R})$-norm decay estimates ($q\in[1,+\infty]$) for two Cauchy systems of type Rao-Nakra sandwich beam with a frictional damping or an infinite memory}
\author[Aissa Guesmia]{Aissa Guesmia}
\maketitle

\pagenumbering{arabic}

\begin{center}
Institut Elie Cartan de Lorraine, UMR 7502, Universit\'e de Lorraine\\
3 Rue Augustin Fresnel, BP 45112, 57073 Metz Cedex 03, France
\end{center}
\begin{abstract}
In this paper, we consider two systems of type Rao-Nakra sandwich beam in the whole line $\mathbb{R}$ with a frictional damping 
or an infinite memory acting on the Euler-Bernoulli equation. When the speeds of propagation of the two wave equations are equal, we show that the solutions do not converge to zero when time goes to infinity. In the reverse situation, we prove some $L^2 (\mathbb{R})$-norm and $L^1 (\mathbb{R})$-norm decay estimates of solutions and theirs higher order derivatives with respect to the space variable. Thanks to interpolation inequalities and Carlson inequality, these $L^2 (\mathbb{R})$-norm and $L^1 (\mathbb{R})$-norm decay estimates lead to similar ones in the $L^q (\mathbb{R})$-norm, for any $q\in [1,+\infty]$. In our both $L^2 (\mathbb{R})$-norm and $L^1 (\mathbb{R})$-norm decay estimates, we specify the decay rates in terms of the regularity of the initial data and the nature of the control.    
\end{abstract}

asymptotic behavior, $L^2 (\mathbb{R})$-norm and 
$L^1 (\mathbb{R})$-norm decay estimates, energy method, Fourier analysis.
\vskip0,1truecm
{\it AMS classification.} 34B05, 34D05, 34H05, 35B40, 35L45, 74H40, 93D20, 93D15.

\renewcommand{\thefootnote}{}
\footnotetext{E-mail addresses: aissa.guesmia@univ-lorraine.fr, ORCID: 0000 0003 4782 5053}

\section{Introduction}

Let $\rho_1 ,\,\rho_2 ,\,\rho_3 ,\,k_0 ,\,k_1,\,k_2 ,\,k_3,\,\gamma$ and $l$ be real positive constants and $g:\mathbb{R}_+:=[0,+\infty) \to \mathbb{R}_+$ satisfying 
$g\in C^1 (\mathbb{R}_{+})$,
\begin{equation}\label{g2}
0<g_0 :=\displaystyle \int_0^{+\infty} g(s) \,ds < k_3 .
\end{equation}
and, for some real positive constants $\beta_1$ and $\beta_2$,
\begin{equation}\label{g1}
-\beta_2 g\leq g^{\prime}\leq -\beta_1 g. 
\end{equation} 
Condition \eqref{g1} implies that $g^{\prime}\leq 0$, $g$ is integrable over 
$\mathbb{R}_+$ and, by integrating,
\begin{equation}
g(0)e^{-\beta_2 s}\leq g(s)\leq g(0)e^{-\beta_1 s}, \quad s\in \mathbb{R}_+ .\label{g3} 
\end{equation}
So \eqref{g2} is valid if $g(0)>0$ and $\beta_1$ is big enough. For example, one can take $g(s)=d_1 e^{-d_2 s}$ with $d_1 ,\,d_2 >0$ satisfying 
$\frac{d_1}{d_2}<k_3$, so \eqref{g2} and \eqref{g1} hold, for $\beta_1 =\beta_2 =d_2$.      
\vskip0,1truecm
This paper deals with the stability of two systems of type Rao-Nakra sandwich beam in the whole line $\mathbb{R}$ with a control acting only on the Euler-Bernoulli equation. These systems consist of two wave equations for the longitudinal displacements of the top and bottom layers, and one Euler-Bernoulli equation for the transversal displacement. The considered control is provided through a frictional damping of size $\gamma$:
\begin{equation}
\left\{
\begin{array}{ll}
\rho_1\varphi_{tt} -k_1 \varphi_{xx} + k_0 (\varphi +\psi+lw_x)=0, \vspace{0.2cm}\\
\rho_2\psi_{tt} -k_2 \psi_{xx} + k_0 (\varphi +\psi+lw_x) =0,
\vspace{0.2cm}\\
\rho_3 w_{tt} +k_3 w_{xxxx} -lk_0 (\varphi +\psi+lw_x)_x + \gamma w_t =0,\vspace{0.2cm}\\
\varphi(x,0) =\varphi_{0} (x),\,\psi(x,0) =\psi_{0} (x),\,w (x,0) =w_{0} (x),\vspace{0.2cm}\\
\varphi_t (x,0) =\varphi_{1} (x),\,\psi_t (x,0) =\psi_{1} (x),\,w_t (x,0) =w_{1} (x)
\end{array}
\right.
\label{S1}
\end{equation}
or an infinite memory of kernel $g$: 
\begin{equation}
\left\{
\begin{array}{ll}
\rho_1\varphi_{tt} -k_1 \varphi_{xx} + k_0 (\varphi +\psi+lw_x)=0, \vspace{0.2cm}\\
\rho_2\psi_{tt} -k_2 \psi_{xx} + k_0 (\varphi +\psi+lw_x)=0,
\vspace{0.2cm}\\
\rho_3 w_{tt} +k_3 w_{xxxx} -lk_0 (\varphi +\psi+lw_x)_x - \displaystyle \int_0^{+\infty} g(s) \, w_{xxxx} (x,t-s) ds =0,
\vspace{0.2cm}\\
\varphi(x,0) =\varphi_{0} (x),\,\psi(x,0) =\psi_{0} (x),\,w (x,-t) =w_{0} (x,t),\vspace{0.2cm}\\
\varphi_t (x,0) =\varphi_{1} (x),\,\psi_t (x,0) =\psi_{1} (x),\,w_t (x,0) =w_{1} (x),
\end{array}
\right.
\label{S2}
\end{equation}
where $x\in \mathbb{R}$ and $t\in \mathbb{R}_+$ are, respectively, the spacial and time variables, $\varphi_{0} ,\,\varphi_{1},\,\psi_0 ,\,\psi_1 ,\,w_0$ and $w_1$ are fixed initial data, and 
\begin{equation*}
(\varphi,\,\psi ,\,w):\,\mathbb{R}\times\mathbb{R}_+\to\mathbb{R}^3 
\end{equation*}
is the unknown of \eqref{S1} and \eqref{S2}. A subscript $r$ denotes the derivative with respect to $r$. Also, we use $\partial_r^m$ or 
$\frac{d^m}{d r^m}$ to denote the differential operator of order $m$ with respect to $r$; i.e. $\frac{\partial^m}{\partial r^m}$. When a function has only one variable, its derivative is noted by $\prime$.
\vskip0,1truecm
Since the works \cite{140, 170, 190}, several layer laminated beam and plate models have been introduced during the last sixty years. The known generalized Rao-Nakra beam (composed of a top and a bottom face plate), presented in \cite{5}, takes into account the shear effect of the bottom and top layers, and it is gien by
\begin{equation}
\left\{
\begin{array}{ll}
\rho_1 h_1 u_{tt} -E_1 h_1 u_{xx} -\tau =0, \vspace{0.2cm}\\
\rho_3 h_3 v_{tt} -E_3 h_3 v_{xx} +\tau =0, \vspace{0.2cm}\\
\rho hw_{tt} +EI w_{xxxx} -G_1 h_1 (w_x +\phi_1 )_x -G_3 h_3 (w_x +\phi_3 )_x - h_2 \tau_x =0,
\vspace{0.2cm}\\
\rho_1 I_1 \phi_{1,tt} -E_1 I_1 \phi_{1,xx} -\frac{1}{2} h_1 \tau +G_1 h_1 (w_x +\phi_1) =0,\vspace{0.2cm}\\
\rho_3 I_3 \phi_{3,tt} -E_3 I_3 \phi_{3,xx} -\frac{1}{2} h_3 \tau +G_3 h_3 (w_x +\phi_3 ) =0,
\end{array}
\right.
\label{GNRB}
\end{equation} 
where $\rho_j ,\,h_j ,\,E_j ,\,G_j ,\,I_j ,\,EI$ and $\rho h$ are real positive constants representing some physical parameters and satisfying some relationships, $u,\,\phi_1 ,\,v$ and $\phi_3$ are, respectively, the longitudinal displacement and shear angle of the top and bottom layers, 
$w$ is the transverse displacement of the beam and $\tau$ is the shear stress in the core layer defined by 
\begin{equation*}
\tau =-u-\frac{1}{2} h_1 \phi_1 +h_2 w_x +v-\frac{1}{2} h_3 \phi_3 . 
\end{equation*}
\vskip0,1truecm
By neglecting some components and/or parameters and/or considering some connections between them, several models were derived from \eqref{GNRB} and studied in the literature like the Rao-Nakra sandwich beam \cite{6}:
\begin{equation}
\left\{
\begin{array}{ll}
\rho_1 h_1 u_{tt} -E_1 h_1 u_{xx} -k(-u+v+\alpha w_x ) =0, \vspace{0.2cm}\\
\rho_3 h_3 v_{tt} -E_3 h_3 v_{xx} +k(-u+v+\alpha w_x ) =0, \vspace{0.2cm}\\
\rho hw_{tt} +EI w_{xxxx} - \alpha k(-u+v+\alpha w_x )_x =0,
\end{array}
\right.
\label{NRB}
\end{equation}
the laminated beam \cite{3}:
\begin{equation}
\left\{
\begin{array}{ll}
\rho w_{tt} +G (\psi -w_{x})_x =0, \vspace{0.2cm}\\
I_{\rho} (3s-\psi)_{tt} -D(3s-\psi)_{xx} -G(\psi -u_x ) =0, \vspace{0.2cm}\\
3I_{\rho} s_{tt} -3Ds_{xx} +3G(\psi -w_x) +4\mu s+4\delta s_t=0,
\end{array}
\right.
\label{LB}
\end{equation}
the Bresse model \cite{bres}:
\begin{equation}
\left\{
\begin{array}{ll} 
\rho_1 u_{tt} -k_1 (u_{x} +\psi +lw)_x -lk_3 (w_x -lu)=0,\vspace{0.2cm}\\
\rho_2 \psi_{tt} -k_2 \psi_{xx} +k_1 (u_x +\psi +lw) =0,\vspace{0.2cm}\\
\rho_3 w_{tt} -k_3 (w_{x} -lu)_x +lk_1 (u_x +\psi +lw) =0,
\end{array}
\right.
\label{BB}
\end{equation}
the Mead-Markos model \cite{140}: 
\begin{equation}
\left\{
\begin{array}{ll} 
\rho hw_{tt} +EI w_{xxxx} -\alpha_1 (-u+v+\alpha w_x ) =0,\vspace{0.2cm}\\
(-u+v+\alpha w_x )_{xx} -\alpha_2 (-u+v+\alpha w_x )-\alpha w_{xxxx} =0
\end{array}
\right.
\label{MMB}
\end{equation} 
and the Timoshenko beam \cite{1}:
\begin{equation}
\left\{
\begin{array}{ll} 
\rho_1 u_{tt} -k_1 (u_{x} +\psi )_x =0,\vspace{0.2cm}\\
\rho_2 \psi_{tt} -k_2 \psi_{xx} +k_1 (u_x +\psi) =0,
\end{array}
\right.
\label{TB}
\end{equation}
where all the coefficients are real postive constants.
\vskip0,1truecm
The study of long time behavior of \eqref{NRB}-\eqref{TB} has been the subject of an active research and mathematical endeavor. Thereby, a huge number of research articles have been appeared. In order to highlight the main contribution and the novelty of the present paper, we
will merely point out the articles whose content concerns the Rao-Nakra sandwich beam \eqref{NRB}, and so their content is closely related to our problems \eqref{S1} and \eqref{S2}. For \eqref{LB}-\eqref{TB} (in different contexts), we refer the reader, for example, to \cite{allen1, allen2, cdflr, gues2, gues3, gues1, 4, hali, 5, 2, s2, 16} and the references therein. 
\vskip0,1truecm
In \cite{7}, the authors considered the Rao-Nakra sandwich beam with 
$x\in (0,1)$ and an internal frictional damping acting either on the beam equation or one of the wave equations, and proved that the polynomial stability occurs. Similar polynomial stability results were proved in \cite{8} under internal damping or Kelvin-Voigt damping working on two of the three equations.  
\vskip0,1truecm
The authors of \cite{rvfp} proved the well-posedness and exponential stability of Rao-Nakra sandwich beam in $(0,L)$, $L>0$, under a heat conduction of secound sound acting on the second wave equation (see 
\cite{17, 18, 19} for details on the second sound mechanism). The same 
well-posedness and exponential stability results of \cite{rvfp} were proved in \cite{13} using boundary feedbacks at $x=L$.   
\vskip0,1truecm
The (global or local) boundary controllability problems of the Rao-Nakra beam were also the subject of several studies in the literature; see, for example, \cite{9, 10, 11, 12, 14, 15} and the references therein. 
\vskip0,1truecm
The subject of this paper is to treat the stability of the Rao-Nakra sandwich beam \eqref{NRB} in the whole line $\mathbb{R}$ under a frictional damping \eqref{S1} or an infinite memory \eqref{S2} (where we denoted $(u,v)$ by $(-\varphi,\psi)$ and simplified the notation of the coefficients). The frictinal damping and infinite memory terms generate the unique dissipation for \eqref{S1} and \eqref{S2} (see \eqref{Ep21} and Remark \ref{remark10} below). This dissipation is acting on the third equation of \eqref{S1} and \eqref{S2} (transversal displacement). To the best of our knowledge, these situations have never been considered in the literature. 
\vskip0,1truecm
The main objective of this paper is to get decay estimates of 
\begin{equation}\label{L2L1Norms}
t\mapsto \left\Vert \partial _{x}^{j}U(\cdot,t)\right\Vert_{2} \quad\hbox{and}\quad t\mapsto \left\Vert \partial _{x}^{j}U(\cdot,t)\right\Vert_{1}, 
\end{equation}
where $U$ is defined in \eqref{UU0} below, $j\in\mathbb{N}$ and 
$\Vert\cdot\Vert_q$ denotes the norm of $L^q (\mathbb{R})$, for $q\in [1,+\infty]$. We will prove the instability of \eqref{S1} and \eqref{S2} when  
\begin{equation}
\frac{k_1}{\rho_1} = \frac{k_2}{\rho_2}. \label{taukrho1}
\end{equation}
However, when
\begin{equation}
\frac{k_1}{\rho_1} \ne\frac{k_2}{\rho_2} , \label{taukrho2}
\end{equation}  
we prove that the functions \eqref{L2L1Norms} satisfy some polynomial stability estimates with decay rates depending on the regularity of the initial data.
\vskip0,1truecm
Our $L^2 (\mathbb{R})$-norm and $L^1 (\mathbb{R})$-norm decay estimates cover all the $L^q (\mathbb{R})$-norm decay estimates of $\partial _{x}^{j}U$, for any $q\in [1,+\infty]$. Indeed, using interpolation inequalities, the decay estimates of the $L^{\infty}$-norm is obtained immediately from the inequality ($\delta_q$ denotes a real positive constant, which does not depend neither on $U$ nor on $j$) 
\begin{equation} \label{Interpol_1}
\Vert \partial _{x}^{j}U\Vert_{\infty}\leq \delta_{\infty} \Vert \partial
_{x}^{j}U\Vert_{2}^{1/2}\Vert \partial _{x}^{j+1}U\Vert_{2}^{1/2}.
\end{equation}
This eventually lead to the decay rate of the $L^q$-norm, 
$2< q< +\infty$, through the interpolation inequality
\begin{equation} \label{Interpol_2}
\Vert \partial _{x}^{j}U\Vert_{q}\leq \delta_q\Vert \partial
_{x}^{j}U\Vert_{\infty}^{1-2/q}\Vert \partial _{x}^{j}U\Vert_{2}^{2/q}.
\end{equation}
To fill the gap $1<q< 2$, we use the interpolation inequality 
\begin{equation}\label{Interpol_3}
\left\Vert \partial _{x}^{j}U\right\Vert _{{q}}\leq \delta_q\left\Vert
\partial _{x}^{j}U\right\Vert _{{2}}^{2(q-1)/q}\left\Vert \partial
_{x}^{j}U\right\Vert _{1}^{(2-q)/q}.
\end{equation}
To get the $L^1 (\mathbb{R})$-norm decay estimate, we treat the asymptotic behavior of $\Vert x\partial_{x}^{j}U\Vert_{{2}}$, and then we use the Carlson inequality (see \cite{BarzBurPe_1998} for instance) 
\begin{equation}
\Vert \partial_{x}^{j}U\Vert_{{1}}\leq \Vert \partial_{x}^{j}U\Vert _{{2}}^{1/2}\Vert x\partial_{x}^{j}U\Vert_{{2}}^{1/2} . \label{Carlson}
\end{equation}     
\vskip0,1truecm
The proof is based on the energy method combined with the Fourier analysis (by using the transformation in the Fourier space).
\vskip0,1truecm
The paper is organized as follows. In section 2, we formulate \eqref{S1} and \eqref{S2} in a first order Cauchy system. Section 3 will be devoted to the proof of the asymptotic behavior of $\vert \widehat{U}\vert$. In section 4, we prove our $L^2 (\mathbb{R})$-norm decay estimate. The asymptotic behavior of $\vert \partial_{\xi}\widehat{U}\vert$ will be treated in section 5. Finally, in section 6, we prove our 
$L^1 (\mathbb{R})$-norm decay estimate. The last section presents some concluding remarks.

\section{Abstract formulation of \eqref{S1} and \eqref{S2}}

To formulate \eqref{S1} and \eqref{S2} in an abstract first order system, we consider the following variables:
\begin{equation}\label{v21}
u=\varphi_t ,\quad y=\psi_t ,\quad \theta =w_t ,\quad v=\varphi_{x} ,\quad z=\psi_x , \quad \phi =w_{xx} \quad\hbox{and}\quad p=\varphi +\psi +lw_x .
\end{equation}
We observe that \eqref{S1}$_1$-\eqref{S1}$_3$ can be written in the form
\begin{align}\label{e10}
\begin{cases}
v_{t} - u_{x} = 0, \\
\rho_1 u_{t} - k_1 v_{x} +k_0 p= 0, \\
z_{t} - y_{x} = 0, \\
\rho_2 y_{t} - k_2 z_{x} + k_0 p = 0, \\
\phi_{t} - \theta_{xx} = 0,  \\
\rho_3\theta_{t} + k_3\phi_{xx} -l k_0 p_x + \gamma \theta = 0,\\
p_t -u-y-l\theta_x =0.
\end{cases}
\end{align}
In case \eqref{S2}, we consider the additional variable introduced in \cite{da}
\begin{align}\label{eta21}
\eta (x,t, s)=w (x,t)-w (x,t-s)
\end{align}
with its initial data $\eta_0 (x,s)=\eta (x,0,s)$. The variable $\eta$ satisfies
\begin{align}\label{eta22}
\eta_t (x,t, s)+\eta_s (x,t, s)=w_t (x,t)
\end{align}
and
\begin{align}\label{eta23}
\displaystyle\int_0^{+\infty}\,g(s) \eta_{xxxx} (x,t,s)\,ds =g_0 w_{xxxx} (x,t)-\displaystyle\int_0^{+\infty}\,g(s)w_{xxxx}(x,t-s)\,ds.
\end{align}
Then \eqref{S2}$_1$-\eqref{S2}$_3$ can be formulated in the form 
\begin{align}\label{e20}
\begin{cases}
v_{t} - u_{x} = 0, \\
\rho_1 u_{t} - k_1 v_{x} +k_0 p= 0, \\
z_{t} - y_{x} = 0, \\
\rho_2 y_{t} - k_2 z_{x} + k_0 p = 0, \\
\phi_{t} - \theta_{xx} = 0,  \\
\rho_3\theta_{t} + (k_3 -\tau_2 g_0)\phi_{xx} -l k_0 p_x + \displaystyle\int_0^{+\infty}\,g(s) \eta_{xxxx} \,ds = 0,\\
p_t -u-y-l\theta_x =0,\\
\eta_t +\eta_s -\theta =0.
\end{cases}
\end{align}
We define the variable $U$ and its initial data $U_0$ by 
\begin{align}\label{UU0}
U (x, t)= \begin{cases}
(v,\,u,\,z,\,y,\,\phi,\,\theta ,\,p )^{T} (x, t) \quad&\hbox{in case \eqref{S1}}, \\
(v,\,u,\,z,\,y,\,\phi,\,\theta ,\,p,\eta)^{T} (x, t)\quad&\hbox{in case \eqref{S2}}
\end{cases}\quad\hbox{and}\quad U_0 (x)= U(x, 0).
\end{align}
Therefore, the systems \eqref{e10} and \eqref{e20} lead to 
\begin{align}\label{e11}
\begin{cases}
U_{t} (x,t)+ A U (x,t)  = 0,\\
U (x, 0)=U_{0} (x),
\end{cases}  
\end{align}
where 
\begin{equation*}
A U=A_4 U_{xxxx} +A_2 U_{xx}+A_1 U_x +A_0 U
\end{equation*}
and the operators $A_j$ are defined in case \eqref{e10} by  
\begin{align*}
\begin{cases}
A_4 U_{xxxx} =(0,0,0,0,0,0,0)^T ,\\ 
A_2 U_{xx} = \left(0,0,0,0,- \theta_{xx} ,\frac{k_3}{\rho_3} \phi_{xx} ,0\right)^T ,\\
A_1 U_x =\left(-u_{x} ,-\frac{k_1}{\rho_1} v_{x} ,- y_{x} ,-\frac{k_2}{\rho_2} z_{x} ,0 ,-\frac{l k_0}{\rho_3} p_x ,-l\theta_x \right)^T ,\\
A_0 U = \left(0,\frac{k_0}{\rho_1} p ,0,\frac{k_0}{\rho_2}p ,0,\frac{\gamma}{\rho_3}\theta ,-u-y \right)^T ,
\end{cases}
\end{align*}
and $A_j$ in case \eqref{e20} are defined by 
\begin{align*}
\begin{cases}
A_4 U_{xxxx} = \left(0,0,0,0,0,\frac{1}{\rho_3} \displaystyle\int_0^{+\infty}\,g(s) \eta_{xxxx}\,ds,0,0\right)^T ,\\
A_2 U_{xx} = \left(0,0,0,0,- \theta_{xx} ,\frac{k_3 -g_0 }{\rho_3}\phi_{xx} ,0,0\right)^T,\\
A_1 U_x = \left(-u_{x} ,-\frac{k_1}{\rho_1}v_{x} ,- y_{x} ,-\frac{k_2}{\rho_2}z_{x} ,0,-\frac{l k_0 }{\rho_3}p_x ,-l\theta_x ,0\right)^T ,\\
A_0 U = \left(0,\frac{k_0}{\rho_1} p ,0 ,\frac{k_0}{\rho_2} p,0 ,0 ,-u-y,\eta_s - \theta\right)^T .
\end{cases}
\end{align*}
\vskip0,1truecm
For a given function $h :\mathbb{R}\to \mathbb{C}$, we use the notations 
$Re\, h$, $Im\, h$, $\bar{h}$ and $\widehat{h}$
to denote, respectively, the real part of $h$, the imaginary part of $h$, the conjugate of $h$ and the Fourier transformation of $h$ given by
\begin{equation}
\widehat{h} (\xi)=\displaystyle\int_{-\infty}^{+\infty}\,e^{-i\xi x} h(x)\,dx,\quad \xi\in \mathbb{R}. \label{hathdef}
\end{equation}
Applying the Fourier transformation (with respect to the space variable $x$) to \eqref{e11}, we obtain the following first-order Cauchy system in the Fourier space:
\begin{align}\label{g6}
\begin{cases}
\widehat{U}_{t} (\xi,\,t)+{\tilde A}(\xi) \widehat{U} (\xi,\,t) = 0, \\
\widehat{U}(\xi,\,0) = \widehat{U}_{0}(\xi),
\end{cases}
\end{align}
where ${\tilde A}(\xi)=\xi^4\,{A}_4 -\xi^2\,{A}_2 + i\,\xi\,{A}_1 + A_0$. The solution of \eqref{g6} is given by
\begin{equation}\label{e14}
\widehat{U}(\xi,\,t) = e^{-{\tilde A}(\xi)\,t}\ \widehat{U}_{0}(\xi).
\end{equation}
Computing the term $e^{-{\tilde A}(\xi)t}$ is a challenging problem, and in many situations, this cannot be done. Consequently, in order to show the
asymptotic behavior of $\widehat{U}$, it suffices to find a non negative function $f(\xi)$ and two positive constants ${\tilde c}$ and $c$ such that, for each $(\xi ,t)\in \mathbb{R}\times \mathbb{R}_+$, 
\begin{equation}  \label{Lambda_estimate}
\left\vert e^{-{\tilde A}(\xi)t}\right\vert\leq {\tilde c}e^{-cf(\xi)t}.
\end{equation}
\vskip0,1truecm
Let $\widehat{E}$ be the energy associated with \eqref{g6} given by
\begin{eqnarray}\label{E21}
\widehat{E}(\xi,t)  &=& \frac{1}{2}\left[k_1 \,\vert\widehat{v}\vert^2 + \rho_1\vert\widehat{u}\vert^{2} + k_2 \,\vert\widehat{z}\vert^{2} + \rho_2\vert\widehat{y}\vert^{2} + (k_3 -{\tau}_0 g_0 )\vert\widehat{\phi}\vert^{2} + \rho_3\vert\widehat{\theta}\vert^{2} +k_0 \vert\widehat{p}\vert^{2} \right](\xi,\,t)\nonumber\\
&&+\frac{{\tau}_0}{2} \xi^4\displaystyle\int_{0}^{+\infty}\,g(s)\vert\widehat{\eta}(\xi,\,t)\vert^{2} \,ds,  
\end{eqnarray}
where 
\begin{eqnarray}\label{tau0}
{\tau}_0  =\begin{cases}
0 \quad&\hbox{in case \eqref{e10}}, \\
1 \quad&\hbox{in case \eqref{e20}}.
\end{cases} 
\end{eqnarray}
\vskip0,1truecm
\begin{lemma}
The energy functional $\widehat{E}$ satisfies, for each 
$(\xi ,t)\in \mathbb{R}\times \mathbb{R}_+$, 
\begin{equation}\label{Ep21}
\frac{d}{dt}\widehat{E}(\xi,\,t) = 
-(1-{\tau}_0)\gamma\,\vert\widehat{\theta}(\xi,\,t)\vert^{2} +
\frac{{\tau}_0}{2} \xi^4\displaystyle\int_{0}^{+\infty}\,g^{\prime}(s)\vert\widehat{\eta}(\xi,\,t)\vert^{2} \,ds.
\end{equation}
\end{lemma}
\vskip0,1truecm
\begin{proof}
The equation \eqref{g6}$_1$ in case \eqref{e10} is equivalent to 
\begin{align}\label{e101}
\begin{cases}
\widehat{v}_{t} - i\xi\widehat{u} = 0, \\
\rho_1 \widehat{u}_{t} - ik_1 \xi \,\widehat{v} +k_0 \widehat{p}= 0, \\
\widehat{z}_{t} - i\xi\widehat{y} = 0, \\
\rho_2\widehat{y}_{t} - ik_2\xi\,\widehat{z} + k_0 \widehat{p} = 0, \\
\widehat{\phi}_{t} +\xi^2\widehat{\theta} = 0,  \\
\rho_3\widehat{\theta}_{t} -k_3\xi^2 \widehat{\phi}- ilk_0\,\xi\widehat{p} + \gamma \widehat{\theta} = 0,\\
\widehat{p}_{t} -\widehat{u} -\widehat{y} -il\xi \widehat{\theta} =0.
\end{cases}
\end{align}
Multiplying the equations in \eqref{e101} by $k_1\bar{\widehat{v}}$, 
$\bar{\widehat{u}}$, $k_2\bar{\widehat{z}}$, $\bar{\widehat{y}}$, $k_3\bar{\widehat{\phi}}$, $\bar{\widehat{\theta}}$ and 
$k_0\bar{\widehat{p}}$, respectively, adding the obtained equations, taking the real part of the resulting expression and using the identity, for $h ,\,d :\mathbb{R}\to \mathbb{C}$ two differentiable functions,
\begin{equation}\label{hd1}
\frac{d}{dt} Re\,\left(h {\bar d}\right) = Re\,\left(h_t {\bar d} +d_t {\bar h}\right),
\end{equation}
we obtain \eqref{Ep21} with ${\tau}_0 =0$. Similarily, \eqref{g6}$_1$ in case \eqref{e20} is reduced to 
\begin{align}\label{e102}
\begin{cases}
\widehat{v}_{t} - i\xi\widehat{u} = 0, \\
\rho_1 \widehat{u}_{t} - ik_1 \xi \,\widehat{v} +k_0 \widehat{p}= 0, \\
\widehat{z}_{t} - i\xi\widehat{y} = 0, \\
\rho_2\widehat{y}_{t} - ik_2 \xi\,\widehat{z} + k_0 \widehat{p} = 0, \\
\widehat{\phi}_{t} +\xi^2\widehat{\theta} = 0,  \\
\rho_3\widehat{\theta}_{t} -(k_3 -g_0 )\xi^2 \widehat{\phi}- ilk_0\,\xi\widehat{p} +  \xi^4 \displaystyle\int_{0}^{+\infty}\,g(s)\widehat{\eta} \,ds= 0,\\
\widehat{p}_{t} -\widehat{u} -\widehat{y} -il\xi \widehat{\theta} =0,\\
\widehat{\eta}_{t} +\widehat{\eta}_s - \widehat{\theta} =0.
\end{cases}
\end{align} 
Multiplying \eqref{e102}$_1$-\eqref{e102}$_7$ by $k_1 \bar{\widehat{v}}$, 
$\bar{\widehat{u}}$, $k_2 \bar{\widehat{z}}$, 
$\bar{\widehat{y}}$, $(k_3 - g_0 )\bar{\widehat{\phi}}$, 
$\bar{\widehat{\theta}}$ and $k_0\bar{\widehat{p}}$, respectively, multiplying \eqref{e102}$_8$ by $\xi^4 g(s)\bar{\widehat{\eta}}$ and integrating on $\mathbb{R}_+$ with respect to $s$, adding all the obtained equations, taking the real part of the resulting expression and using \eqref{hd1}, we get 
\begin{equation*}
\frac{d}{dt}\widehat{E}(\xi,\,t) = -\frac{1}{2} \xi^4\displaystyle\int_{0}^{+\infty}\,g(s)\frac{d}{ds}\vert\widehat{\eta}\vert^{2} \,ds,
\end{equation*}
therefore, by integrating with respect to $s$, we obtain   
\begin{equation*}
\frac{d}{dt}\widehat{E}(\xi,\,t) = -\frac{1}{2} \left[\xi^4 \displaystyle\int_{0}^{+\infty}\,g(s)\vert\widehat{\eta}\vert^{2} \,ds\right]_{s=0}^{s=+\infty} +\frac{1}{2} \xi^4\displaystyle\int_{0}^{+\infty}\,g^{\prime}(s)\vert\widehat{\eta}\vert^{2} \,ds.
\end{equation*}
Because \eqref{g3} and \eqref{eta21} imply that  
\begin{equation}\label{limgeta}
\displaystyle\lim_{s\to +\infty}g(s)=0\quad\hbox{and}\quad\widehat{\eta} (\xi,t,0)=0,
\end{equation}
then we find \eqref{Ep21} with ${\tau}_0 =1$.
\end{proof} 
\vskip0,1truecm
\begin{remark}\label{remark10}
Notice that \eqref{Ep21} implies that 
\begin{equation*}
\frac{d}{dt}\widehat{E}(\xi,\,t) \leq 0,
\end{equation*}
since $\gamma>0$ and $g^{\prime} \leq 0$, so \eqref{g6} is dissipative. 
If the frictional damping and infinite memory are not considered
(i.e. $\gamma =g = 0$), then \eqref{g6} is conservative; that is 
\begin{equation*}
\widehat{E}(\xi,\,t) = \widehat{E}(\xi,\,0).
\end{equation*}
On the other hand, we put
\begin{eqnarray*}
\vert\widehat{U}(\xi,\,t)\vert^2  =\left[\vert\widehat{v}\vert^2 + \vert\widehat{u}\vert^{2} + \vert\widehat{z}\vert^{2} +\vert\widehat{y}\vert^{2} + \vert\widehat{\phi}\vert^{2} + \vert\widehat{\theta}\vert^{2} +\vert\widehat{p}\vert^{2} \right](\xi,\,t) +\tau_0 \xi^4\displaystyle\int_{0}^{+\infty}\,g(s)\vert\widehat{\eta}(\xi,\,t)\vert^{2} \,ds. 
\end{eqnarray*}
So we deduce that, thanks to the right inequality in \eqref{g2}, 
\begin{equation}\label{equiv}
\vert\widehat{U}\vert^2\sim \widehat{E} . 
\end{equation} 
\end{remark}

\section{Estimation of $\vert\widehat{U}\vert$}

This section is dedicated to the study of the asymptotic behavior of 
$\widehat{U}(\xi,\,t)$, when time $t$ goes to infinity. We prove the next theorem.
\vskip0,1truecm
\begin{theorem}\label{theorem31} 
Let $\widehat{U}$ be the solution of \eqref{g6}. 
\vskip0,1truecm
1. In case \eqref{taukrho1} and for every $\xi\in \mathbb{R}$, $\vert\widehat{U}(\xi,\,t)\vert$ doesn't converge to zero when $t$ goes to infinity.
\vskip0,1truecm
2. In case \eqref{taukrho2}, there exist $c,\,\widetilde{c}>0$ (independent on $\xi$ and $t$) such that 
\begin{equation}\label{11}
\vert\widehat{U}(\xi,\,t)\vert \leq\widetilde{c}\,e^{-\,c\,f(\xi)\,t}\,\vert\widehat{U}_0 (\xi)\vert,\quad\forall (\xi ,t)\in \mathbb{R}\times \mathbb{R}_+ , 
\end{equation}
where  
\begin{equation}\label{22}
f(\xi) = \frac{\xi^{4+2\tau_0}}{\xi^{10}+1}. 
\end{equation}
\end{theorem}
\vskip0,1truecm

\subsection{Case \eqref{taukrho1}} We prove here, under \eqref{taukrho1} and for every $\xi\in \mathbb{R}$, that $\vert\widehat{U}(\xi,\,t)\vert$ doesn't converge to zero when time $t$ goes to infinity. It suffices to prove that, for any $\xi\in \mathbb{R}$, $-{\tilde A}(\xi)$ has at least a pure imaginary eigenvalue (see \cite{tesc}); that is
\begin{equation}\label{Axi}
\forall \xi\in \mathbb{R},\,\,\exists \lambda\in \mathbb{R}^*,\,\,\exists \widehat{U}\ne 0:\quad i\lambda \widehat{U} + {\tilde A}(\xi)\widehat{U}=0. 
\end{equation}
\vskip0,1truecm
{\bf System \eqref{e101}}. To get \eqref{Axi}, it is enough to prove that 
\begin{equation}\label{Axidet}
\forall \xi\in \mathbb{R},\,\,\exists \lambda\in \mathbb{R}^* :\quad det\left(i\lambda I+{\tilde A}(\xi)\right)=0, 
\end{equation}
where $I$ denotes the identity matrix. We see that \eqref{taukrho1} implies, in case \eqref{e101}, that
\begin{equation*}
i\lambda I+{\tilde A}(\xi) = 
\begin{pmatrix}
i\lambda & -i\xi & 0 & 0 & 0 & 0 & 0 \\
-i\frac{k_1}{\rho_1} \xi & i\lambda & 0 & 0 & 0 & 0 & \frac{k_0}{\rho_1}\\
0 & 0 & i\lambda & -i\xi & 0 & 0 & 0  \\
0 & 0 & -i\frac{k_1}{\rho_1}\xi & i\lambda & 0 & 0 & \frac{k_0}{\rho_2} \\
0 & 0 & 0 & 0 & i\lambda & \xi^2 & 0 \\
0 & 0 & 0 & 0 & -\frac{k_3}{\rho_3} \xi^2 & i\lambda +\frac{\gamma}{\rho_3} & -i\frac{lk_0}{\rho_3}\xi \\
0 & -1 & 0 & -1 & 0 &-il\xi & i\lambda
\end{pmatrix} .
\end{equation*}
A direct computation shows that 
\begin{eqnarray}
det\left(i\lambda I+{\tilde A}(\xi)\right) &=& i\lambda\left(\lambda^2 -\frac{k_1}{\rho_1} \xi^2 \right)\left[\lambda^2 -k_0 \left(\frac{1}{\rho_1}+\frac{1}{\rho_2}\right)-\frac{k_1}{\rho_1} \xi^2\right]\left[i\lambda \left(i\lambda + \frac{\gamma}{\rho_3}\right)+\frac{k_3}{\rho_3}\xi^4 \right]\nonumber\\
&& +i\frac{l^2 k_0}{\rho_3}\lambda \xi^2 \left(\lambda^2 -\frac{k_1}{\rho_1} \xi^2 \right)^2 . \label{detiLIA} 
\end{eqnarray}
It is clear that $det\left(i\lambda I+{\tilde A}(\xi)\right)=0$ for 
\begin{align}\label{lambdaxikjrhoj}
\lambda =\begin{cases}
{\sqrt{\frac{k_1}{\rho_1}}} \xi \quad\quad &\hbox{if}\,\,\xi \ne 0,\\
{\sqrt{k_0 \left(\frac{1}{\rho_1}+\frac{1}{\rho_2}\right)}}\quad\quad &\hbox{if}\,\,\xi = 0. 
\end{cases}
\end{align}
Consequently, \eqref{Axidet} holds.
\vskip0,1truecm
{\bf System \eqref{e102}}. For $\lambda\in \mathbb{R}^*$, we put 
\begin{equation}\label{hatetahattheta}
\widehat{\eta} (s)=\frac{1}{i\lambda}\left(1-e^{-i\lambda s} \right)\widehat{\theta}\quad\hbox{and}\quad {\tilde g} (\lambda)=\int_0^{+\infty}
\left(1-e^{-i\lambda s} \right)g(s)ds. 
\end{equation} 
We observe that ${\tilde g} (\lambda)$ is well defined (according to \eqref{g2}) and $\widehat{\eta}$ is the unique function satisfying \eqref{Axi}$_8$ and $\widehat{\eta} (0)=0$. On the other hand, \eqref{taukrho1} implies that the first seven equations of \eqref{Axi} are equivalent to ${\tilde B}(\xi)\left(\widehat{v}, \widehat{u}, \widehat{z},\widehat{y},\widehat{\phi},\widehat{p}\right)^T =0$, where
\begin{equation*}
{\tilde B}(\xi) = \begin{pmatrix}
i\lambda & -i\xi & 0 & 0 & 0 & 0 & 0 \\
-i\frac{k_1}{\rho_1} \xi & i\lambda & 0 & 0 & 0 & 0 & \frac{k_0}{\rho_1}\\
0 & 0 & i\lambda & -i\xi & 0 & 0 & 0  \\
0 & 0 & -i\frac{k_1}{\rho_1}\xi & i\lambda & 0 & 0 & \frac{k_0}{\rho_2} \\
0 & 0 & 0 & 0 &  i\lambda & \xi^2 & 0 \\
0 & 0 & 0 & 0 & -\frac{k_3 -g_0}{\rho_3} \xi^2 & i\lambda +\frac{{\tilde g} (\lambda)}{\rho_3}\xi^4 & -i\frac{lk_0}{\rho_3}\xi \\
0 & -1 & 0 & -1 & 0 &-il\xi & i\lambda
\end{pmatrix} .
\end{equation*} 
Then the problem \eqref{Axi} is reduced to prove that
\begin{equation}\label{detB}
\forall \xi\in \mathbb{R},\,\,\exists \lambda\in \mathbb{R}^* :\quad det\,{\tilde B}(\xi)=0. 
\end{equation}
A direct computation shows that (as for \eqref{detiLIA} with $k_3 -g_0$ and 
$\xi^4 {\tilde g} (\lambda)$ instead of $k_3$ and $\gamma$, respectively) 
\begin{eqnarray*}
det {\tilde B}(\xi) &=& i\lambda\left(\lambda^2 -\frac{k_1}{\rho_1} \xi^2 \right)\left[\lambda^2 -k_0 \left(\frac{1}{\rho_1}+\frac{1}{\rho_2}\right)-\frac{k_1}{\rho_1} \xi^2\right]\left[i\lambda \left(i\lambda + \frac{{\tilde g} (\lambda)}{\rho_3}\xi^4\right)+\frac{k_3 -g_0 }{\rho_3}\xi^4 \right]\nonumber\\
&& +i\frac{l^2 k_0}{\rho_3}\lambda \xi^2 \left(\lambda^2 -\frac{k_1}{\rho_1} \xi^2 \right)^2 .  
\end{eqnarray*} 
We remark that $det\,{\tilde B}(\xi)=0$ for $\lambda$ given by \eqref{lambdaxikjrhoj}. Thus \eqref{detB} holds. Consequently, the proof of the first result in Theorem \ref{theorem31} is ended. 

\subsection{Some differential equations} We start the proof of \eqref{11} by proving some useful differential equations. To simplify the computations, we put 
\begin{equation}\label{G12}
G = (1-\tau_0 )\gamma \widehat{\theta}+\tau_0 \xi^4 \int_0^{+\infty} g(s)\widehat{\eta}ds, 
\end{equation}
where $\tau_0$ is defined in \eqref{tau0}. Let consider the system
\begin{align}\label{SG12}
\begin{cases}
\widehat{v}_{t} - i\xi\widehat{u} = 0, \\
\rho_1 \widehat{u}_{t} - ik_1 \xi \,\widehat{v} +k_0 \widehat{p}= 0, \\
\widehat{z}_{t} - i\xi\widehat{y} = 0, \\
\rho_2\widehat{y}_{t} - ik_2 \xi\,\widehat{z} + k_0 \widehat{p} = 0, \\
\widehat{\phi}_{t} +\xi^2\widehat{\theta} = 0,  \\
\rho_3\widehat{\theta}_{t} -(k_3 -\tau_0g_0 )\xi^2 \widehat{\phi}- ilk_0\,\xi\widehat{p} + G = 0,\\
\widehat{p}_{t} -\widehat{u} -\widehat{y} -il\xi \widehat{\theta} =0,\\
\widehat{\eta}_{t} +\widehat{\eta}_s - \widehat{\theta} =0.
\end{cases}
\end{align}
It is evident that \eqref{e101} is identical to \eqref{SG12}$_1$-\eqref{SG12}$_7$ if $\tau_0 =0$, and \eqref{e102} coincides with \eqref{SG12} if 
$\tau_0 =1$. 
\vskip0,1truecm
Multiplying \eqref{SG12}$_{4}$ and \eqref{SG12}$_{3}$ by $i\,\xi\,\overline{\widehat{z}}$ and $-i\rho_2\,\xi\,\overline{\widehat{y}}$, 
respectively, adding the resulting equations, taking the real part and using \eqref{hd1}, we obtain
\begin{equation}\label{eq31}
\frac{d}{dt}Re\left(i\rho_2\,\xi\,\widehat{y}\, \overline{\widehat{z}}\right) = \xi^{2}\left(\rho_2\vert\widehat{y}\vert^{2} -k_2 \vert\widehat{z}\vert^{2}\right) + k_0\,Re\left(i\,\xi\,\widehat{z}\,\overline{\widehat{p}}\right) .
\end{equation}
Similarily, multiplying \eqref{SG12}$_{2}$ and \eqref{SG12}$_{1}$ by $i\,\xi\,\overline{\widehat{v}}$ and $-i\rho_1\,\xi\,\overline{\widehat{u}}$, 
respectively, adding the resulting equations, taking the real part and using \eqref{hd1}, we find
\begin{equation}\label{eq32}
\frac{d}{dt}Re\left(i\rho_1\,\xi\,\widehat{u}\, \overline{\widehat{v}}\right) = \xi^{2}\left(\rho_1 \vert\widehat{u}\vert^{2} -k_1 \vert\widehat{v}\vert^{2}\right) + k_0 Re\left(i\,\xi\,\widehat{v}\,\overline{\widehat{p}}\right) .
\end{equation}
Also, multiplying \eqref{SG12}$_{6}$ and \eqref{SG12}$_{5}$ by 
$-\overline{\widehat{\phi}}$ and $-\rho_3\,\overline{\widehat{\theta}}$, respectively, adding the resulting equations, taking the real part and using \eqref{hd1}, we get
\begin{equation}\label{eq33}
\frac{d}{dt}Re\left(-\rho_3\,\widehat{\theta}\, \overline{\widehat{\phi}}\right) = \xi^{2}\left(\rho_3\vert\widehat{\theta}\vert^{2} -(k_3 -\tau_0 g_0 ) \vert\widehat{\phi}\vert^{2}  \right) +lk_0\, Re\left(i\,\xi\,\widehat{\phi}\,\overline{\widehat{p}}\right) +\,Re\left( \overline{\widehat{\phi}}\,G\right).
\end{equation}
Multiplying \eqref{SG12}$_{4}$ and \eqref{SG12}$_{7}$ by $\xi^2\,\overline{\widehat{p}}$ and $\rho_2\xi^2\,\overline{\widehat{y}}$, 
respectively, adding the resulting equations, taking the real part and using \eqref{hd1}, we obtain
\begin{equation}\label{Sandeq35}
\frac{d}{dt}Re\left(\rho_2\xi^2\,\widehat{p}\, \overline{\widehat{y}}\right) = \xi^{2}\left(\rho_2\,\vert\widehat{y}\vert^{2} - k_0\vert\widehat{p}\vert^{2}\right) + \rho_2 \xi^{2}\,Re\left(\widehat{y}\,\overline{\widehat{u}}\right)+ l\rho_2 \xi^{2}\,Re\left(i\xi\,\widehat{\theta}\,\overline{\widehat{y}}\right)  + k_2\xi^{2}\,Re\left(i\,\xi\,\widehat{z}\,\overline{\widehat{p}}\right).
\end{equation}
After, multiplying \eqref{SG12}$_{3}$ and \eqref{SG12}$_{6}$ by $\rho_3\,\overline{\widehat{\theta}}$ and $\overline{\widehat{z}}$, respectively, adding the resulting equations, taking the real part and using \eqref{hd1}, we entail
\begin{equation}\label{eq36}
\frac{d}{dt}Re\left(\rho_3\widehat{z}\, \overline{\widehat{\theta}}\right) = -\rho_3\, Re\left(i\xi\widehat{\theta}\,\overline{\widehat{y}}\right)+(k_3 -\tau_0g_0)\,\xi^{2}\, Re\left(\widehat{\phi}\,\overline{\widehat{z}}\right)-lk_0\,Re\left(i\,\xi\,\widehat{z}\,\overline{\widehat{p}}\right)
-\,Re\left(\overline{\widehat{z}}\, G\right).
\end{equation}
Multiplying \eqref{SG12}$_{5}$ and \eqref{SG12}$_{4}$ by $i\rho_2 \xi\,\overline{\widehat{y}}$ and $-i \xi\overline{\widehat{\phi}}$, respectively, adding the resulting equations, taking the real part and using \eqref{hd1}, it follows that
\begin{equation}\label{eq37}
\frac{d}{dt}Re\left(i\rho_2 \xi\,\widehat{\phi}\, \overline{\widehat{y}}\right) = -\rho_2\xi^{2}\, Re\left(i \xi\widehat{\theta}\,\overline{\widehat{y}}\right)+k_2 \xi^2\, Re\left(\widehat{\phi}\,\overline{\widehat{z}}\right)-k_0\,Re\left(i \xi\widehat{\phi}\,\overline{\widehat{p}}\right).
\end{equation}
Next, multiplying \eqref{SG12}$_{2}$ and \eqref{SG12}$_{3}$ by $i \xi\overline{\widehat{z}}$ and $-i\rho_1 \xi\,\overline{\widehat{u}}$, 
respectively, adding the resulting equations, taking the real part and using \eqref{hd1}, it appears that
\begin{equation}\label{eq310}
\frac{d}{dt}Re\left(i\rho_1 \xi\,\widehat{u}\, \overline{\widehat{z}}\right) = -k_1 \xi^2\, Re\left(\widehat{v}\,\overline{\widehat{z}}\right) +\rho_1  \xi^2\,Re\left(\widehat{y} \,\overline{\widehat{u}}\right)+k_0 \,Re\left(i \xi\widehat{z}\, \overline{\widehat{p}}\right).
\end{equation}
Multiplying \eqref{SG12}$_{2}$ and \eqref{SG12}$_{5}$ by $i \xi\overline{\widehat{\phi}}$ and $-i\rho_1 \xi\,\overline{\widehat{u}}$, respectively, adding the resulting equations, taking the real part and using \eqref{hd1}, we see that
\begin{equation}\label{eq312}
\frac{d}{dt}Re\left(i\rho_1 \xi\,\widehat{u}\, \overline{\widehat{\phi}}\right) = -k_1 \xi^2\,Re\left(\widehat{\phi}\,\overline{\widehat{v}}\right) +\rho_1 \,\xi^2\,Re\left(i \xi\widehat{\theta} \,\overline{\widehat{u}}\right)+k_0\,Re\left(i \xi\widehat{\phi}\, \overline{\widehat{p}}\right).
\end{equation}
Also, multiplying \eqref{SG12}$_{1}$ and \eqref{SG12}$_{4}$ by $i\rho_2 \xi\,\overline{\widehat{y}}$ and $-i \xi\overline{\widehat{v}}$, respectively, adding the resulting equations, taking the real part and using \eqref{hd1}, we find
\begin{equation}\label{eq36}
\frac{d}{dt}Re\left(i\rho_2 \xi\widehat{v}\, \overline{\widehat{y}}\right) =- \rho_2 \xi^2\, Re\left(\widehat{y}\,\overline{\widehat{u}}\right)+k_2 \xi^2\, Re\left(\widehat{v}\,\overline{\widehat{z}}\right)-k_0\,Re\left(i \xi\widehat{v}\,\overline{\widehat{p}}\right).
\end{equation}
Multiplying \eqref{SG12}$_{1}$ and \eqref{SG12}$_{6}$ by $\rho_3\,\overline{\widehat{\theta}}$ and $\overline{\widehat{v}}$, respectively, adding the resulting equations, taking the real part and using \eqref{hd1}, we get
\begin{equation}\label{eq37}
\frac{d}{dt}Re\left(\rho_3\,\widehat{v}\, \overline{\widehat{\theta}}\right) = -\rho_3\, Re\left(i\xi\widehat{\theta}\,\overline{\widehat{u}}\right)+(k_3 -\tau_0 g_0 )\, \xi^2\,Re\left(\widehat{\phi}\,\overline{\widehat{v}}\right)-lk_0\,Re\left(i\xi\widehat{v}\,\overline{\widehat{p}}\right)
-\,Re\left(\overline{\widehat{v}}\, G\right).
\end{equation} 
Now, multiplying \eqref{SG12}$_{6}$ by $-\xi^2 \displaystyle\int_{0}^{+\infty}\,g(s)\overline{\widehat{\eta}} \,ds$, multiplying \eqref{SG12}$_{8}$ by $-\rho_3 \xi^2 g(s)\overline{\widehat{\theta}}$ and integrating over $\mathbb{R}_+$ with respect to $s$, adding the resulting equations, taking the real part and using \eqref{hd1}, we infer that 
\begin{eqnarray}\label{eq6341}
\frac{d}{dt}Re\left(-\rho_3 \xi^2 {\widehat{\theta}}\, \displaystyle\int_{0}^{+\infty}\,g(s)\overline{\widehat{\eta}} \,ds\right) &=& -\rho_3 g_0 \xi^2 \vert \widehat{\theta}\vert^2  -(k_3 -\tau_0 g_0)\xi^4 Re\left(\widehat{\phi}\displaystyle\int_{0}^{+\infty}\,g(s)\overline{\widehat{\eta}} \,ds\right) \nonumber \\ 
&&-lk_0 \xi^2 Re\left(i\xi\widehat{p}\displaystyle\int_{0}^{+\infty}\,g(s)\overline{\widehat{\eta}} \,ds\right) -\rho_3 \xi^2  Re\left(\overline{\widehat{\theta}}\displaystyle\int_{0}^{+\infty}\,g^{\prime}(s){\widehat{\eta}} \,ds\right)\\
&&+\xi^2 Re \left(G\displaystyle\int_{0}^{+\infty}\,g(s)\overline{{\widehat{\eta}}} \,ds\right). \nonumber 
\end{eqnarray}
Similarily, multiplying \eqref{SG12}$_{4}$ by 
$i\xi \displaystyle\int_{0}^{+\infty}\,g(s)\overline{\widehat{\eta}} \,ds$, multiplying \eqref{SG12}$_{8}$ by $-i\rho_2 \xi g(s)\,\overline{\widehat{y}}$ and integrating over $\mathbb{R}_+$ with respect to $s$, adding the resulting equations, taking the real part and using \eqref{hd1}, we get
\begin{eqnarray}\label{Sandeq634}
\frac{d}{dt}Re\left(i\rho_2 \xi {\widehat{y}}\, \displaystyle\int_{0}^{+\infty}\,g(s)\overline{\widehat{\eta}} \,ds\right) &=& -\rho_2 g_0 Re\left(i\xi{\widehat{\theta}} \overline{\widehat{y}}\right) -k_2 \xi^2 Re\left(\widehat{z}\displaystyle\int_{0}^{+\infty}\,g(s)\overline{\widehat{\eta}} \,ds\right) \nonumber \\ 
&&-k_0 Re\left(i\xi\widehat{p}\displaystyle\int_{0}^{+\infty}\,g(s)\overline{\widehat{\eta}} \,ds\right) -\rho_2 Re\left(i\xi\overline{\widehat{y}}\displaystyle\int_{0}^{+\infty}\,g^{\prime}(s){\widehat{\eta}} \,ds\right).
\end{eqnarray}
Finally, multiplying \eqref{SG12}$_{2}$ by $i \xi\displaystyle\int_{0}^{+\infty}\,g(s)\overline{\widehat{\eta}} \,ds$, multiplying \eqref{SG12}$_{8}$ by $-\rho_1 \xi g(s)\overline{\widehat{u}}$ and integrating over $\mathbb{R}_+$ with respect to $s$, adding the resulting equations, taking the real part and using \eqref{hd1}, we see that
\begin{eqnarray}\label{eq6342}
\frac{d}{dt}Re\left(i\rho_1 \xi {\widehat{u}}\, \displaystyle\int_{0}^{+\infty}\,g(s)\overline{{\widehat{\eta}}}\,ds\right) &=& -k_1 \xi^2 Re\left(\widehat{v}\displaystyle\int_{0}^{+\infty}\,g(s)\overline{\widehat{\eta}} \,ds\right)-k_0 Re\left(i \xi\widehat{p}\displaystyle\int_{0}^{+\infty}\,g(s)\overline{\widehat{\eta}} \,ds\right) \nonumber \\ 
&& -\rho_1 Re\left(i \xi\overline{\widehat{u}}\displaystyle\int_{0}^{+\infty}\,g^{\prime}(s){\widehat{\eta}} \,ds\right)-\rho_1 g_0 Re\left(i \xi{\widehat{\theta}}\overline{\widehat{u}}\right) .
\end{eqnarray}   
\vskip0,1truecm
Let $\lambda_1 ,\cdots,\lambda_{13}$ be real numbers to be chosen later (which do not depend on time $t$ but may depend on $\xi$ and the parameters of \eqref{S1} and \eqref{S2}). We define the functionals ${F}_0 ,\,F_1 ,\,F_2 $ and ${F}_3$ as follows:
\begin{eqnarray}\label{F0}
{F}_0 (\xi,\,t) & = & Re\left(i\rho_2 \lambda_1\xi\widehat{y}\overline{\widehat{z}}+i\rho_1 \lambda_2\xi\widehat{u}\overline{\widehat{v}}
-\rho_3 \lambda_3\widehat{\theta}\overline{\widehat{\phi}}+\rho_2 \lambda_4\xi^2\widehat{p}\overline{\widehat{y}}+\rho_3 \lambda_5\widehat{z}\overline{\widehat{\theta}}\right) \\
& & +Re\left(i\rho_2 \lambda_6\xi\widehat{\phi}\overline{\widehat{y}}
+i\rho_1 \lambda_7\xi\widehat{u}\overline{\widehat{z}}+i\rho_1 \lambda_{8}\xi\widehat{u}\overline{\widehat{\phi}}+i\rho_2\lambda_{9}\xi\widehat{v}\overline{\widehat{y}}+\rho_3 \lambda_{10}\widehat{v}\overline{\widehat{\theta}}\right) \nonumber \\
& & +\tau_0 Re\left[\left(-\rho_3 \lambda_{11}\xi^2\widehat{\theta}+i\rho_2 \lambda_{12}\xi\widehat{y}+i\rho_1 \lambda_{13}\xi\widehat{u}\right)\displaystyle\int_{0}^{+\infty}\,g(s)\overline{\widehat{\eta}}ds\right] , \nonumber 
\end{eqnarray}
\begin{equation}\label{F1}
{F}_1 (\xi,\,t) = -\xi^2 \left(B_1 \vert\widehat{z}\vert^{2}+B_2 \vert\widehat{v}\vert^{2} +B_3 \vert\widehat{\phi}\vert^{2} +B_4 \vert\widehat{u}\vert^{2} +B_5 \vert\widehat{p}\vert^{2}+B_6 \vert\widehat{y}\vert^{2}+B_7 \vert\widehat{\theta}\vert^{2}\right) , 
\end{equation}
\begin{equation}\label{F2}
{F}_2 (\xi,\,t) = Re\left[i\xi \left(A_1\widehat{z}+A_2\widehat{v}+A_3 \widehat{\phi}\right)\overline{\widehat{p}}+A_4\widehat{\phi}\overline{\widehat{z}}+A_5\widehat{v}\overline{\widehat{z}}+A_6\widehat{\phi}\overline{\widehat{v}}+A_7\xi^2\widehat{y}\overline{\widehat{u}}+iA_8\xi\widehat{\theta}\overline{\widehat{u}}+iA_9\xi\widehat{\theta}\overline{\widehat{y}}\right] 
\end{equation}
and
\begin{eqnarray}\label{F3}
{F}_3 (\xi,\,t) &=& Re\left[\left(\lambda_3\overline{\widehat{\phi}} -\lambda_5\overline{\widehat{z}}-\lambda_{10}\overline{\widehat{v}}\right)G\right] -\tau_0 Re\left[\left(\rho_3 \lambda_{11}\xi^2 \overline{\widehat{\theta}}+i\rho_2\lambda_{12} \xi \overline{\widehat{y}}+i\rho_1 \lambda_{13}\xi \overline{\widehat{u}}\right)\displaystyle\int_{0}^{+\infty}\,g^{\prime}(s){\widehat{\eta}}ds\right]\nonumber \\
& & +\tau_0 Re\left[\lambda_{11}\left[-(k_3 -\tau_0 g_0 ) \xi^4 \widehat{\phi}-ilk_0 \xi^3\widehat{p}+\xi^2 G\right] \displaystyle\int_{0}^{+\infty}\,g(s)\overline{\widehat{\eta}}ds\right] \nonumber \\
& & -\tau_0 Re\left[\left[\lambda_{12}\left(k_2 \xi^2 \widehat{z}+ik_0 \xi\widehat{p} \right) +\lambda_{13}\left(k_1 \xi^2 \widehat{v}+ik_0 \xi\widehat{p}\right)\right]\displaystyle\int_{0}^{+\infty}\,g(s)\overline{\widehat{\eta}}ds\right],   
\end{eqnarray}
where 
\begin{equation*}
B_1 =k_2\lambda_1 ,\quad B_2 =k_1\lambda_2 ,\quad B_3 =(k_3 -\tau_0 g_0)\lambda_3 ,\quad B_4 =-\rho_1 \lambda_2 ,
\end{equation*}
\begin{equation*}
B_5 =k_0 \lambda_4 ,\quad B_6 =-\rho_2 (\lambda_1 +\lambda_4 ) ,\quad B_7 =\rho_3 (\tau_0 g_0\lambda_{11} -\lambda_3 ),
\end{equation*}
\begin{equation*}
A_1 =k_0 \lambda_1 +k_2 \lambda_{4}\xi^2 -lk_0 \lambda_5 +k_0\lambda_7 ,\quad A_2 =k_0 (\lambda_2 -\lambda_{9} -l\lambda_{10}), \quad A_3 =k_0 (l\lambda_3 -\lambda_{6} +\lambda_{8} ) 
\end{equation*}
\begin{equation*}
A_4 =(k_3 -\tau_0 g_0 )\lambda_{5}\xi^2 +k_2 \lambda_6 \xi^2 ,\quad 
A_5 =-k_1 \lambda_7 \xi^2 +k_2 \lambda_{9}\xi^2, \quad A_6 =-k_1 \lambda_{8} \xi^2 +(k_3 -\tau_0 g_0 )  \lambda_{10}\xi^2 ,
\end{equation*}
\begin{equation*}
A_7 =\rho_1 \lambda_{7}+\rho_2 (\lambda_{4} -\lambda_{9}) ,\quad A_8 =\rho_1 (\lambda_{8}\xi^2 -\tau_0 g_0 \lambda_{13})-\rho_3 \lambda_{10} \quad 
\hbox{and}\quad A_9 =\rho_2 (l\lambda_{4} -\lambda_6 )\xi^2 -\rho_3 \lambda_5 -\tau_0 \rho_2 g_0 \lambda_{12}.
\end{equation*}
Multiplying \eqref{eq31}-\eqref{eq6342} by $\lambda_{1},\cdots,\lambda_{10} ,\,\tau_0 \lambda_{11},\,\tau_0 \lambda_{12}$ and 
$\tau_0\lambda_{13}$, respectively, and adding the obtained equations, we deduce that
\begin{equation} \label{g26+}
\frac{d}{dt}{F}_0 (\xi,\,t)={F}_1 (\xi,\,t) +{F}_2 (\xi,\,t)+{F}_3 (\xi,\,t).
\end{equation}
According to the notations \eqref{v21}, we have $p_x =v+z+l\phi$. 
Because \eqref{hathdef} implies that ${\widehat{p_x}}=i\xi {\widehat{p}}$, then  
\begin{equation*}
i\xi \overline{\widehat{p}}=-\overline{\widehat{v}}-\overline{\widehat{z}}-l\overline{\widehat{\phi}},
\end{equation*}
this identity allows to formulate the first term in $F_2$ as 
\begin{eqnarray}\label{ixihatbarp}
Re\left[i\xi \left(A_1\widehat{z}+A_2\widehat{v}+A_3 \widehat{\phi}\right)\overline{\widehat{p}}\right] &=& -\left(A_1 \vert\widehat{z}\vert^{2}+A_2 \vert\widehat{v}\vert^{2} +lA_3 \vert\widehat{\phi}\vert^{2}\right)-\left(A_1 +A_2 \right)Re\left( {\widehat{v}}\overline{\widehat{z}}\right)\\
& &-\left(A_3 +lA_1 \right)Re \left({\widehat{\phi}}\overline{\widehat{z}}\right)-\left(lA_2 +A_3 \right)Re \left({\widehat{\phi}}\overline{\widehat{v}}\right) \nonumber.
\end{eqnarray}
By combining \eqref{g26+} and \eqref{ixihatbarp}, we get
\begin{equation}\label{dtF1}
\frac{d}{dt}{F}_0 (\xi,\,t)={F}_5 (\xi,\,t) +{F}_4 (\xi,\,t)+{F}_3 (\xi,\,t), 
\end{equation}
where 
\begin{eqnarray*}
{F}_5 (\xi,\,t) & = & -\left[\left(B_1\xi^2 +A_1\right) \vert\widehat{z}\vert^{2}+\left(B_2\xi^2 +A_2\right) \vert\widehat{v}\vert^{2} +\left(B_3 \xi^2 +lA_3\right) \vert\widehat{\phi}\vert^{2} \right] \\
& & -\xi^2 \left[B_4 \vert\widehat{u}\vert^{2} +B_5 \vert\widehat{p}\vert^{2}+B_6 \vert\widehat{y}\vert^{2}+B_7 \vert\widehat{\theta}\vert^{2}\right] 
\end{eqnarray*}
and 
\begin{eqnarray*}
{F}_4 (\xi,\,t) & = & Re\left[\left(A_4 -A_3 -lA_1\right)\widehat{\phi}\overline{\widehat{z}}+\left(A_5 -A_1 -A_2\right)\widehat{v}\overline{\widehat{z}}+\left(A_6 -lA_2 -A_3 \right)\widehat{\phi}\overline{\widehat{v}}\right]\\
& & +Re\left[A_7\xi^2\widehat{y}\overline{\widehat{u}}+iA_8\xi\widehat{\theta}\overline{\widehat{u}}+iA_9\xi\widehat{\theta}\overline{\widehat{y}}\right]. 
\end{eqnarray*}
Now, we choose the different real numbers $\lambda_{1},\cdots,\lambda_{15}$ in order to have
\begin{equation}\label{A1A6}
A_4 -A_3 -lA_1 =A_5 -A_1 -A_2 =A_6 -lA_2 -A_3 =A_7 = 0, 
\end{equation}
\begin{equation}\label{B1B3}
B_1\xi^2 +A_1 ={\tilde B}_1 \xi^2,\quad B_2\xi^2 +A_2 ={\tilde B}_2 \xi^2 ,\quad B_3 \xi^2 +lA_3 ={\tilde B}_3 \xi^2 ,  
\end{equation} 
\begin{equation}\label{tildeB1tildeB3}
{\tilde B}_1 >0 ,\quad {\tilde B}_2 >0,\quad {\tilde B}_3 >0,\quad B_4 >0,\quad B_5 >0,\quad B_6 >0 
\end{equation} 
and
\begin{align}\label{B6B7}
A_8 =A_9 =0\,\,\,\hbox{and}\,\,\, B_7 >0\quad \hbox{if}\,\,\tau_0 =1.
\end{align} 
We start by choosing $\lambda_7 ,\,\lambda_{8}$ and $\lambda_{9}$ as follows:
\begin{equation}\label{lambda91011}
\lambda_7 =k_2 \xi^2 +l\lambda_5 -\lambda_1  ,\quad \lambda_{8} =-\frac{2(k_3 -\tau_0 g_0)}{l}\xi^2 +\lambda_6 -l\lambda_3 \quad\hbox{and}\quad \lambda_{9} =-k_1 \xi^2 +\lambda_2 -l\lambda_{10} . 
\end{equation}  
The choices \eqref{lambda91011} guarantee \eqref{B1B3} with 
\begin{equation}\label{tildeB1tildeB3=}
{\tilde B}_1 =B_1 +k_2 \lambda_4 +k_0 k_2 ,\quad {\tilde B}_2 =B_2 +k_0 k_1 \quad\hbox{and}\quad {\tilde B}_3 =B_3 -2k_0 (k_3 -\tau_0 g_0).
\end{equation}
Also, we select $\lambda_1 ,\,\lambda_5 ,\,\lambda_{6}$ and $\lambda_{10}$  by 
\begin{equation}\label{lambda8}
\lambda_6 =\frac{1}{k_2}\left[lk_2\lambda_4 -(k_3 -\tau_0 g_0 )\lambda_5 +k_0 \left(lk_2 -\frac{2(k_3 -\tau_0 g_0)}{l}\right)\right],
\end{equation}  
\begin{equation}\label{lambda12}
\lambda_{10} =\frac{1}{lk_2}\left[-2k_1 k_2\xi^2 +k_1\lambda_1 +k_2\lambda_2 -k_2 \lambda_4 - lk_1 \lambda_5 -k_0 (k_1 +k_2)\right],
\end{equation}
\begin{equation}\label{lambda6}
\lambda_1 = -\frac{k_2}{k_1}\lambda_2 
-\frac{l^2 k_2 }{k_3 -\tau_0 g_0 }\lambda_3 +\left(\frac{l^2 k_2}{k_3 -\tau_0 g_0}+\frac{k_2}{k_1} \right)\lambda_4 +k_0 \left(\frac{2l^2 k_2}{k_3 -\tau_0 g_0}-1 -\frac{k_2}{k_1}\right)
\end{equation} 
and
\begin{equation}\label{lambda71}
\lambda_{5} =-\frac{k_2}{l}\xi^2 +\frac{1}{l}\lambda_1 -\frac{k_0\rho_2 (k_1 +k_2 )}{l(k_1\rho_2 -k_2\rho_1)}
\end{equation}
($\lambda_{5}$ and $\lambda_{1}$ are well defined thanks to \eqref{taukrho2} and the right inequality in \eqref{g2}, respectively). According to \eqref{lambda91011}, the selections \eqref{lambda8} and \eqref{lambda12} imply $A_4 -A_3 -lA_1 =0$ and $A_5 -A_1 -A_2 =0$, respectively, \eqref{lambda8}-\eqref{lambda6} lead to $A_6 -lA_2 -A_3 =0$, and \eqref{lambda12} and \eqref{lambda71} guarantee 
$A_7 =0$, so \eqref{A1A6} is satisfied. We put
\begin{equation}\label{B0}
B_0 =\min\left\{{\tilde B}_1 ,{\tilde B}_2 ,{\tilde B}_3 ,B_4 ,B_5 ,B_6\right\}.
\end{equation}
Then, multiplying \eqref{dtF1} by $\xi^2$ and exploiting \eqref{A1A6} and \eqref{B1B3}, we obtain
\begin{eqnarray}\label{xiF0}
\frac{d}{dt}\left[\xi^2 {F}_0 (\xi,\,t)\right] & \leq & -B_0\xi^4\left[\vert\widehat{z}\vert^{2}+\vert\widehat{v}\vert^{2} +\vert\widehat{\phi}\vert^{2} +\vert\widehat{u}\vert^{2} +\vert\widehat{p}\vert^{2} +\vert\widehat{y}\vert^{2}\right]-B_7\xi^4\vert\widehat{\theta}\vert^{2} \\
& & +\xi^2 Re\left[iA_8\xi\widehat{\theta}\overline{\widehat{u}}+iA_9\xi\widehat{\theta}\overline{\widehat{y}}\right] +\xi^2 F_3 (\xi ,t).\nonumber 
\end{eqnarray}
\vskip0,1truecm
In order to estimate $\xi^2 F_3 (\xi ,t)$, first, we notice the following evident inequality: 
\begin{equation}\label{m123}
\vert\xi\vert^{m_2}\leq \vert\xi\vert^{m_1} +\vert\xi\vert^{m_3},\quad 
\forall \xi\in\mathbb{R},\,\,\forall 0\leq m_1 \leq m_2\leq m_3 .
\end{equation} 
Second, using H\"older's inequality and the right inequality in \eqref{g1}, we see that, for ${\tilde\eta}\in\left\{{\widehat{\eta}},\overline{\widehat{\eta}}\right\}$, 
\begin{eqnarray*}
\left\vert\displaystyle\int_{0}^{+\infty}\,g(s){\tilde\eta} (\xi ,s)\,ds\right\vert^2 & = & \left\vert\displaystyle\int_{0}^{+\infty}\,{\sqrt{g(s)}}{\sqrt{g(s)}}{\tilde\eta}(\xi ,s) \,ds\right\vert^2 \\
& \leq & \left(\displaystyle\int_{0}^{+\infty}\,g(s)\,ds\right)
\displaystyle\int_{0}^{+\infty}\,g(s)\vert{\tilde\eta}(\xi ,s)\vert^2 \,ds \\
& \leq & -\frac{g_0}{\beta_1} \displaystyle\int_{0}^{+\infty}\,g^{\prime}(s)\vert{\tilde{\eta}}(\xi ,s)\vert^2 \,ds.
\end{eqnarray*} 
Similarly, using the limit in \eqref{limgeta}, we have
\begin{eqnarray*}
\left\vert\displaystyle\int_{0}^{+\infty}\,g^{\prime}(s){\tilde\eta} (\xi ,s)\,ds\right\vert^2 
& = & \left\vert\displaystyle\int_{0}^{+\infty}\,{\sqrt{-g^{\prime}(s)}}{\sqrt{-g^{\prime}(s)}}{\tilde\eta} (\xi ,s) \,ds\right\vert^2 \\
& \leq & \left(-\displaystyle\int_{0}^{+\infty}\,g^{\prime} (s)\,ds\right)
\displaystyle\int_{0}^{+\infty}\,(-g^{\prime}(s))\vert{\tilde\eta}(\xi ,s)\vert^2 \,ds \\
& \leq & -g(0)\displaystyle\int_{0}^{+\infty}\,g^{\prime}(s)\vert{\tilde{\eta}} (\xi ,s)\vert^2 \,ds.
\end{eqnarray*} 
Then, using these two inequalities and Young's inequality, we get, for any 
$h :\mathbb{R}\to \mathbb{C}$ and $\varepsilon >0$, 
\begin{equation}\label{hd41}
Re\,\left(h(\xi)\displaystyle\int_{0}^{+\infty}\,g(s){\tilde\eta} (\xi ,s)\,ds\right)\leq \varepsilon \vert h(\xi)\vert^2 - \frac{g_0}{4\varepsilon\beta_1} \displaystyle\int_{0}^{+\infty}\,g^{\prime}(s)\vert {\tilde{\eta}}(\xi ,s)\vert^{2} \,ds
\end{equation}
and 
\begin{equation}\label{hd42}
Re\,\left(h(\xi)\displaystyle\int_{0}^{+\infty}\,g^{\prime} (s){\tilde\eta} (\xi, s)\,ds\right)\leq \varepsilon \vert h(\xi)\vert^2 - \frac{g(0)}{4\varepsilon}\displaystyle\int_{0}^{+\infty}\,g^{\prime}(s)\vert {\tilde{\eta}}(\xi, s)\vert^{2} \,ds.
\end{equation}
Third, in the sequel, we use $C$ (sometimes $C_1 ,\,C_2 ,\,\cdots$) to denote a generic real positive constant, and $C_{\varepsilon}$ to denote a generic real positive constant depending on some real positive constant 
$\varepsilon$, where $C$ and $C_{\varepsilon}$ may be different from step to step. The constants $C$ ($C_1 ,\,C_2 ,\,\cdots$), $\varepsilon$ and $C_{\varepsilon}$ are independent on $x$, $\xi$ and $t$.
\vskip0,1truecm
Applying Young's inequality for all the terms in $\xi^2 {F}_3 (\xi,\,t)$ and using \eqref{m123}, \eqref{hd41}, \eqref{hd42} and the fact that 
$\tau_0^2 =\tau_0$ and $\tau_0 (1-\tau_0) =0$, we get, for any 
$\varepsilon >0$,
\begin{eqnarray}\label{xiF3}
\xi^2 {F}_3 (\xi,\,t) & \leq & \varepsilon \xi^4 \left[\vert\widehat{z}\vert^{2}+\vert\widehat{v}\vert^{2} +\vert\widehat{\phi}\vert^{2} + \tau_0\left(\vert\widehat{u}\vert^{2} +\vert\widehat{y}\vert^{2}+\vert\widehat{\theta}\vert^{2}+\vert\widehat{p}\vert^{2}\right)\right]\nonumber\\
& & +C_{\varepsilon}\left(\lambda_3^2 +\lambda_5^2 +\lambda_{10}^2 \right)\vert{G}\vert^{2} +\tau_0\vert\lambda_{11}\vert \xi^4\vert G\vert\left\vert\displaystyle\int_{0}^{+\infty}\,g(s)\overline{\widehat{\eta}} (\xi ,s)\,ds\right\vert\\
& & -\tau_0 C_{\varepsilon}\left[\lambda_{11}^2 \xi^4\left(\xi^4 +1\right)+\left(\lambda_{12}^2 +\lambda_{13}^2\right) \xi^2\left(\xi^2 +1\right)\right]\displaystyle\int_{0}^{+\infty}\,g^{\prime}(s)\vert{\widehat{\eta}} (\xi ,s)\vert^2 \,ds \nonumber
\end{eqnarray} 
\vskip0,1truecm      
Now, we select $\lambda_2 ,\,\lambda_3 ,\,\lambda_4 ,\,\lambda_{11},\,\lambda_{12}$ and $\lambda_{13}$ in order to get  \eqref{tildeB1tildeB3} and \eqref{B6B7} (and then, in particular, 
$B_0 >0$ since \eqref{B0}). To do so, we distinguish the two cases related 
to \eqref{tau0}. 

\subsection{Case 1: $\tau_0 =0$} Notice that, because $\tau_0 =0$, the numbers $\lambda_{11},\,\lambda_{12}$ and $\lambda_{13}$ are not used. On the other hand, we select $\lambda_2 ,\,\lambda_3$ and $\lambda_4$ as follows:
\begin{equation}\label{lambda11}
-k_0 <\lambda_2 < 0 ,
\end{equation}
\begin{equation}\label{lambda110000}
0 <\lambda_{4} < \frac{lk_2}{k_1 k_3 +k_2\left(l^2 k_1 +k_3\right)}\left[\frac{k_3}{l}\lambda_2 +lk_1 \lambda_3 + k_0\left(\frac{k_1 k_3}{lk_2} +\frac{k_3}{l}-2l k_1 \right)\right],
\end{equation}
\begin{equation}\label{lambda110000+}
\lambda_{4} > \frac{k_1 k_3}{k_1 k_3 +k_2\left(l^2 k_1 +k_3\right)}\left[\frac{k_2}{k_1}\lambda_2 +\frac{l^2 k_2}{k_3}\lambda_3 +k_0\left(\frac{k_2}{k_1} -\frac{2l^2 k_2}{k_3}\right)\right]
\end{equation}
and
\begin{equation}\label{lambda610000}
\lambda_3 > \max\left\{2k_0 ,\frac{1}{l^2 k_1 k_2}\left[-k_2 k_3\lambda_2 +k_0 \left(2l^2 k_1 k_2 -k_1 k_3 -k_2 k_3\right)\right]\right\}. 
\end{equation}
We go back in the reverse order to select the numbers in the following manner: 
\vskip0,1truecm 
1. Choose $\lambda_2$ by \eqref{lambda11}.   
\vskip0,1truecm 
2. After, select $\lambda_3$ large enough in such a way 
that \eqref{lambda610000} holds. 
\vskip0,1truecm 
3. Take $\lambda_4$ so that \eqref{lambda110000} and \eqref{lambda110000+} are true. According to \eqref{lambda610000}, $\lambda_4$ exists.   
\vskip0,1truecm 
4. Now, it is possible to find $\lambda_1$ through \eqref{lambda6}. 
\vskip0,1truecm
5. After, take $\lambda_5$ as in \eqref{lambda71}. 
\vskip0,1truecm
6. Next, it is time to pick $\lambda_6$ and $\lambda_{10}$ verifying \eqref{lambda8} and \eqref{lambda12}, respectively.
\vskip0,1truecm
7. Finally, we can select $\lambda_7$, $\lambda_{8}$ and $\lambda_{9}$ by 
\eqref{lambda91011}.  
\\
\\
We observe that the left inequality in \eqref{lambda11}, \eqref{lambda610000}, the right inequality in \eqref{lambda11} and the left one in \eqref{lambda110000} imply that ${\tilde B}_2 >0,\,{\tilde B}_3 >0,\,B_4 >0$ and $B_5 >0$, respectively. Moreover, \eqref{lambda110000+} and the right inequality in \eqref{lambda110000} combined with \eqref{lambda6} imply 
${\tilde B}_1 >0$ and $B_6 >0$, respectively. Therefore \eqref{tildeB1tildeB3} holds. Thus \eqref{A1A6}-\eqref{tildeB1tildeB3} are satisfied. 
\vskip0,1truecm
On the other hand, because $\lambda_1 ,\,\lambda_2 ,\,\lambda_3$ and $\lambda_4$ do not depend on $\xi$, 
\begin{equation}\label{vertlambdajleq}
\vert\lambda_j\vert\leq C\left(\xi^2 +1\right),\,\, j=5,\cdots,10,\quad\hbox{and}\quad \vert A_j\vert\leq C\left(\xi^4 +1\right),\,\, j=8,9
\end{equation}
(since \eqref{m123}), then, applying Young's inequality and using \eqref{m123}, we get, for any $\varepsilon >0$, 
\begin{equation}\label{xiF021}
\xi^2 Re\left[iA_8\xi\widehat{\theta}\overline{\widehat{u}}+iA_9\xi\widehat{\theta}\overline{\widehat{y}}\right]\leq\varepsilon\xi^4\left(\vert\widehat{u}\vert^{2} +\vert\widehat{y}\vert^{2}\right)+C_{\varepsilon}\xi^2\left(\xi^8 +1\right)\vert\widehat{\theta}\vert^{2} .
\end{equation}
Moreover, we conclude from \eqref{G12} and \eqref{xiF3} (with $\tau_0 =0$) that   
\begin{equation}\label{xiF031}
\xi^2 {F}_3 (\xi,\,t) \leq \varepsilon \xi^4 \left[\vert\widehat{z}\vert^{2}+\vert\widehat{v}\vert^{2} +\vert\widehat{\phi}\vert^{2}\right] +C_{\varepsilon}\left(\xi^4 +1\right)\vert\widehat{\theta}\vert^{2} .
\end{equation}
By combining \eqref{xiF0}, \eqref{xiF021} and \eqref{xiF031}, choosing 
$0<\varepsilon <B_0$ ($\varepsilon$ exists thanks to \eqref{tildeB1tildeB3} and \eqref{B0}), and using \eqref{E21} (with $\tau_0 =0$) and \eqref{m123}, we deduce that there exist real positive constants $c_1$ and $c_2$ such that (notice that $B_7$ does not depend on $\xi$) 
\begin{equation}\label{xiF04100}
\frac{d}{dt}\left[\xi^2 {F}_0 (\xi,\,t)\right] \leq -c_1 \xi^4 \widehat{E}(\xi,\,t)+c_2 \left(\xi^{10} +1\right)\vert\widehat{\theta}\vert^{2}. 
\end{equation} 
Now, let $\lambda$ be a real positive constant and 
\begin{equation}\label{xiL1100}
F(\xi,\,t) = \lambda\,\widehat{E}(\xi,\,t) + \frac{\xi^2}{\xi^{10} +1}\,{F}_0 (\xi,\,t) .
\end{equation}
From \eqref{Ep21} (with $\tau_0 =0$), \eqref{xiF04100} and \eqref{xiL1100}, we find 
\begin{equation}\label{xiL2100}
\frac{d}{dt}{F}(\xi,\,t) \leq -c_1 f(\xi)\widehat{E}(\xi,\,t)-\left(\gamma\,\lambda - c_2\right)\vert\widehat{\theta}\vert^{2} , 
\end{equation}
where $f$ is defined by \eqref{22} with $\tau_0 =0$. Moreover, using the definitions of $\widehat{E} ,\,{F}_0$ and ${F}$, and exploiting \eqref{m123} and \eqref{vertlambdajleq}, we observe that there exists $c_{3} >0$ (not depend on $\lambda$) such that 
\begin{equation}\label{xiL3} 
\vert{F}(\xi,\,t) - \lambda\,\widehat{E}(\xi,\,t)\vert\leq \frac{c_{3}\,\xi^2 (\vert\xi\vert^3 +1 )}{\xi^{10} +1}\,\widehat{E}(\xi,\,t) \leq 2c_{3}\,\widehat{E}(\xi,\,t).
\end{equation} 
Therefore, for $\lambda$ satisfying $\lambda >\max\left\{\frac{c_2}{\gamma},\,2c_{3}\right\}$, we deduce from \eqref{xiL2100} and \eqref{xiL3} that 
\begin{equation}\label{xiL4}
\frac{d}{dt}{F}(\xi,\,t) \leq - c_1 \,f(\xi)\,\widehat{E}(\xi,\,t) 
\end{equation}
and $F\sim \widehat{E}$, since
\begin{equation}\label{xiL5}
(\lambda -2c_{3})\,\widehat{E}(\xi,\,t) \leq {F}(\xi,\,t) \leq (\lambda +2c_{3})\,\widehat{E}(\xi,\,t).
\end{equation}
Consequently, a combination of \eqref{xiL4} and the second inequality in \eqref{xiL5} leads to, for $c=\frac{c_1}{2(\lambda +2c_{3})}$,  
\begin{equation}\label{xiL6}
\frac{d}{dt}{F} (\xi,\,t) \leq -2c\,f(\xi)\,F (\xi,\,t) .
\end{equation}
Multiplying \eqref{xiL6} by $e^{2cf(\xi)t}$, we find 
\begin{equation}\label{xiL6L6}
\frac{d}{dt} \left(e^{2cf(\xi)t}{F} (\xi,\,t)\right)\leq 0,
\end{equation}
then, by integration \eqref{xiL6L6} with respect to $t$,
\begin{equation}\label{xiL6L6L6}
{F} (\xi,\,t)\leq e^{-2cf(\xi)t}{F} (\xi,\,0).
\end{equation}
Finally, using \eqref{equiv} and \eqref{xiL5}, \eqref{11} in case 
$\tau_0 =0$ follows from \eqref{xiL6L6L6}. 

\subsection{Case 2: $\tau_0 =1$} We choose $\lambda_1 ,\cdots,\,\lambda_{10}$ as in case $1$ but with $k_3 -g_0$ instead of $k_3$ ($k_3 -g_0 >0$ thanks to the right inequality in \eqref{g2}), and we get \eqref{A1A6}-\eqref{tildeB1tildeB3}. Next, we select
\begin{equation}\label{lambda1504}
\lambda_{13} =\frac{1}{g_0}\lambda_{8} \xi^2 -\frac{\rho_3}{g_0 \rho_1}\lambda_{10},
\end{equation}
\begin{equation}\label{lambda1404}
\lambda_{12} =\frac{1}{g_0}(l\lambda_{4} -\lambda_6)\xi^2 -\frac{\rho_3}{g_0\rho_2}\lambda_{5} 
\end{equation}
and
\begin{equation}\label{lambda1304}
\lambda_{11} > \frac{1}{g_0}\lambda_{3} 
\end{equation}
($g_0 >0$ according to the left inequality in \eqref{g2}). The choices \eqref{lambda1504}-\eqref{lambda1304} imply tha $A_8 =A_9 =0$ and $B_7 >0$, respectively, thus \eqref{B6B7} is satisfied. Therefore, \eqref{xiF0} implies 
\begin{equation}\label{xiF0A8A9B7}
\frac{d}{dt}\left[\xi^2 {F}_0 (\xi,\,t)\right] \leq  -\min\{B_0 ,B_7\}\xi^4\left[\vert\widehat{z}\vert^{2}+\vert\widehat{v}\vert^{2} +\vert\widehat{\phi}\vert^{2} +\vert\widehat{u}\vert^{2} +\vert\widehat{p}\vert^{2} +\vert\widehat{y}\vert^{2}+\vert\widehat{\theta}\vert^{2}\right]+\xi^2 F_3 (\xi ,t).
\end{equation}
On the other hand, $\lambda_1 ,\,\lambda_2 ,\,\lambda_3 ,\lambda_4$ and $\lambda_{11}$ do not depend on $\xi$, then 
\begin{equation}\label{vertlambdajleq+}
\vert\lambda_j\vert\leq C\left(\xi^2 +1\right),\,\, j=5,\cdots,10,\quad\hbox{and}\quad \vert \lambda_j\vert\leq C\left(\xi^4 +1\right),\,\, j=12,13,
\end{equation}
therefore, applying Young's inequality and using \eqref{G12} (with 
$\tau_0 =1$) and \eqref{m123}, we conclude from \eqref{xiF3}, for any $\varepsilon >0$, that   
\begin{eqnarray}\label{xiF0324}
\xi^2 {F}_3 (\xi,\,t) &\leq& \varepsilon \xi^4 \left[\vert\widehat{z}\vert^{2}+\vert\widehat{v}\vert^{2} +\vert\widehat{\phi}\vert^{2} +\vert\widehat{p}\vert^{2}+\vert\widehat{u}\vert^{2}+\vert\widehat{y}\vert^{2}+\vert\widehat{\theta}\vert^{2}\right] \nonumber\\ 
&& -C_{\varepsilon}\xi^2\left(\xi^{10} +1\right)\displaystyle\int_{0}^{+\infty}\,g^{\prime}(s)\vert{\widehat{\eta}} (\xi ,s)\vert^2 \,ds.
\end{eqnarray}
By combining \eqref{xiF0A8A9B7} and \eqref{xiF0324}, choosing 
$0<\varepsilon <\min\{B_0 ,B_7\}$, and using \eqref{E21} and \eqref{m123}, we deduce that there exist real positive constants $c_1$ and $c_2$ such that 
\begin{equation}\label{xiF0414}
\frac{d}{dt}\left[\xi^2 {F}_0 (\xi,\,t)\right] \leq -c_1 \xi^4 \widehat{E}(\xi,\,t)-c_2 \xi^2 \left(\xi^{10} +1\right)\displaystyle\int_{0}^{+\infty}\,g^{\prime}(s)\vert{\widehat{\eta}} (\xi ,s)\vert^2 \,ds. 
\end{equation} 
Now, let $\lambda$ be a real positive constant and 
\begin{equation}\label{xiL1144}
F(\xi,\,t) = \lambda\,\widehat{E}(\xi,\,t) + \frac{\xi^4}{\xi^{10} +1}\,{F}_0 (\xi,\,t) .
\end{equation}
From \eqref{Ep21} (with $\tau_0 =1$), \eqref{xiF0414} and \eqref{xiL1144}, we find 
\begin{equation}\label{xiL210004}
\frac{d}{dt}{F}(\xi,\,t) \leq -c_1 f(\xi)\widehat{E}(\xi,\,t)+\left(\frac{\lambda}{2} - c_2\right)\xi^4\displaystyle\int_{0}^{+\infty}\,g^{\prime}(s)\vert{\widehat{\eta}} (\xi ,s)\vert^2 \,ds , 
\end{equation}
where $f$ is defined in \eqref{22} with $\tau_0 =1$. Also, using the definitions of $\widehat{E} ,\,{F}_0$ and ${F}$, and exploiting \eqref{m123} and \eqref{vertlambdajleq+}, we see that there exists $c_{3} >0$ (not depend on $\lambda$) such that 
\begin{equation}\label{xiL3000} 
\vert{F}(\xi,\,t) - \lambda\,\widehat{E}(\xi,\,t)\vert\leq \frac{c_{3}\,\xi^4\left(\vert\xi\vert^5 +1\right)}{\xi^{10} +1}\,\widehat{E}(\xi,\,t) \leq 2c_{3}\,\widehat{E}(\xi,\,t).
\end{equation} 
Therefore, choosing $\lambda >\max\left\{2c_2 ,\,2c_{3}\right\}$ and using \eqref{xiL210004} and \eqref{xiL3000} and the fact that $g^{\prime}\leq 0$, we find \eqref{xiL4} and \eqref{xiL5}. Consequently, \eqref{xiL4} and the second inequality in \eqref{xiL5} lead to \eqref{xiL6}.
Finally, by integration \eqref{xiL6} with respect to $t$ and using \eqref{equiv} and \eqref{xiL5}, we obtain \eqref{11} in case $\tau_0 =1$. The proof of the second result in Theorem \ref{theorem31} is finished.

\section{Estimation of $\left\Vert\partial_{x}^{j}U\right\Vert_{2}$}

In this section, we use the second result in Theorem \ref{theorem31} to get some decay estimates on $\Vert\partial_{x}^{k}U\Vert_{2}$ when \eqref{taukrho2} is valid. 
\vskip0,1truecm
\begin{theorem}\label{theorem41}
Assume that \eqref{taukrho2} is satisfied. Let $N\,\in \mathbb{N}^*$, 
\begin{equation}\label{U0regularity}
U_{0}\in H^{N}(\mathbb{R})\cap L^{1}(\mathbb{R})
\end{equation}
and $U$ be the solution of \eqref{e11}. Then, for any $\ell\in\{1,\,2,\,\ldots,\,N\}$ and $j\in\{0,\,1,\,\ldots,\,N-\ell\}$, there exists $c_0 >0$ such that 
\begin{equation}\label{sandwi332}
\Vert\partial_{x}^{j}U\Vert_{2} \leq c_0 \,(1 + t)^{-1/8 - j/4}\,\Vert U_{0} \Vert_{1} + c_0 \,(1 + t)^{-\ell/6}\,\Vert\partial_{x}^{j+\ell}U_{0} \Vert_{2},\quad t\in \mathbb{R}_+ 
\end{equation}
if $\tau_0 =0$, and 
\begin{equation}\label{sandwi334}
\Vert\partial_{x}^{j}U\Vert_{2} \leq c_0 \,(1 + t)^{-1/12 - j/6}\,\Vert U_{0} \Vert_{1} + c_0 \,(1 + t)^{-\ell/4}\,\Vert\partial_{x}^{j+\ell}U_{0} \Vert_{2},\quad t\in \mathbb{R}_+ 
\end{equation}
if $\tau_0 =1$.  
\end{theorem}
\vskip0,1truecm
\begin{proof} For the proof of \eqref{sandwi332}, we see that \eqref{22} with $\tau_0 =0$ implies that (low and high frequences)
\begin{equation}\label{Sg36}
f(\xi) \geq \left\{
\begin{array}{cc}
\frac{1}{2}\,\xi^{4} & \mbox{if} \quad \vert\xi \vert\leq 1, \\
\\
\frac{1}{2}\,\xi^{-6} & \mbox{if} \quad \vert\xi\vert > 1.
\end{array}
\right.
\end{equation}
Applying Plancherel's theorem and \eqref{11}, we entail
\begin{equation}\label{Sg37}
\Vert\partial_{x}^{j}U\Vert_{2}^{2} = \left\Vert\widehat{\partial_{x}^{j} U}(x,\,t)\right \Vert_{2}^{2}  = \int_{\mathbb{R}}\xi^{2\,j}\,\vert\widehat{U}(\xi,\,t)\vert^{2} d\xi 
\end{equation}
\begin{eqnarray*}
& \leq & \widetilde{c}\int_{\mathbb{R}}\xi^{2\,j}\,e^{-\,c\,f(\xi)\,t}\,\vert\widehat{U}_0 (\xi)\vert^{2} d\xi \\
& \leq & \widetilde{c}\int_{\vert\xi\vert\leq 1}\xi^{2\,j}\,e^{-\,c\,f(\xi)\,t}\,\vert\widehat{U}_0 (\xi)\vert^{2}  d\xi + \widetilde{c}\int_{\vert\xi\vert > 1}\xi^{2\,j}\,e^{-\,c\,f(\xi)\,t}\,\vert\widehat{U}_0 (\xi)\vert^{2}\, d\xi := J_{1} + J_{2}.
\end{eqnarray*}
Applying Lemma 2.3 of \cite{gues2} and \eqref{Sg36}$_1$, 
it follows, for the low frequency region,
\begin{equation}\label{Sg38}
J_{1} \leq C\,\Vert\widehat{U}_{0}\Vert_{\infty}^{2}\int_{\vert\xi\vert\leq 1}\xi^{2\,j}\,e^{-\,\frac{c}{2}\,t\,\xi^{4}}\ d\xi \leq C\left(1 + t\right)^{-\,\frac{1}{4}(1 + 2\,j)}\Vert {U}_{0}\Vert_{1}^{2}.
\end{equation}
For the high frequency region, using \eqref{Sg36}$_2$, we observe that
\begin{eqnarray*}
J_{2} & \leq & C\int_{\vert\xi\vert > 1}\vert\xi\vert^{2\,j}\,e^{-\frac{c}{2}\,t\,\xi^{-6}}\,\vert\widehat{U}(\xi,\,0)\vert^{2}\, d\xi \\
& \leq & C\ \sup_{\vert\xi\vert >1}\left\{\vert\xi\vert^{-2\,\ell}\, e^{-\frac{c}{2}\,t\,\vert\xi\vert^{-6}}\right\} \int_{\mathbb{R}}\vert\xi\vert^{2\,(j + \ell)}\,\vert\widehat{U}(\xi,\,0)\vert^{2}\ d\xi,
\end{eqnarray*}
then, using Lemma 2.4 of \cite{gues2},
\begin{equation}\label{Sg38+}
J_{2} \leq C\left(1 + t\right)^{-\frac{1}{3}\,\ell}\,\Vert\,\partial_{x}^{j + \ell}U_{0}\,\Vert_{2}^{2} , 
\end{equation}
and so, by combining \eqref{Sg37}-\eqref{Sg38+}, we get 
\begin{equation}\label{Sg38+000}
\Vert\partial_{x}^{j}U\Vert_{2}^{2} \leq C\left[\left(1 + t\right)^{-\,\frac{1}{4}(1 + 2\,j)}\Vert {U}_{0}\Vert_{1}^{2} +\left(1 + t\right)^{-\frac{1}{3}\,\ell}\,\Vert\,\partial_{x}^{j + \ell}U_{0}\,\Vert_{2}^{2}  \right]. 
\end{equation}
Finally, by combining \eqref{Sg38+000} and the inequality
\begin{equation}\label{sqrta1a2}
{\sqrt{a_1 +a_2}}\leq {\sqrt{a_1}} +{\sqrt{a_2}},\quad\forall a_1 ,a_2\in \mathbb{R}_+ ,
\end{equation}
we find \eqref{sandwi332}.
\vskip0,1truecm
The proof of \eqref{sandwi334} is identical to the one of \eqref{sandwi332} where, instead of \eqref{Sg36}, we have (according to \eqref{22} with 
$\tau_0 =1$)
\begin{equation}\label{Sg360}
f(\xi) \geq \left\{
\begin{array}{cc}
\frac{1}{2}\,\xi^{6} & \mbox{if} \quad \vert\xi \vert\leq 1, \\
\\
\frac{1}{2}\,\xi^{-4} & \mbox{if} \quad \vert\xi\vert > 1.
\end{array}
\right.
\end{equation}  
\end{proof}

\section{Estimation of $\vert \partial_{\xi}\widehat{U}\vert$}

In this section, we study the asymptotic behavior (with respect to $t$) of 
$\partial_{\xi}\widehat{U}$. In order to simplify the computations, let us denoting $\partial_{\xi}\widehat{U}=\widehat{\mathbf U}$, $\partial_{\xi}\widehat{U}_0=\widehat{\mathbf U}_0$ and
\begin{equation*}
\left(\partial_{\xi}\widehat{v},\partial_{\xi}\widehat{u},\partial_{\xi}\widehat{z},\partial_{\xi}\widehat{y},\partial_{\xi}\widehat{\phi},\partial_{\xi}\widehat{\theta},\partial_{\xi}\widehat{p},\partial_{\xi}\widehat{\eta}\right)=\left(\widehat{\mathbf v},\widehat{\mathbf u},\widehat{\mathbf z},\widehat{\mathbf y},\widehat{\Phi},\widehat{\Theta},\widehat{\mathbf p},\widehat{\Lambda}\right).  
\end{equation*}
As \eqref{E21}, the energy associated to $\widehat{\mathbf U}$ is define by
\begin{eqnarray}\label{E21part}
\widehat{\mathbf E}(\xi,t)  &=& \frac{1}{2}\left[k_1 \,\vert\widehat{\mathbf v}\vert^2 + \rho_1\vert\widehat{\mathbf u}\vert^{2} + k_2 \,\vert\widehat{\mathbf z}\vert^{2} + \rho_2\vert\widehat{\mathbf y}\vert^{2} + (k_3 -{\tau}_0 g_0 )\vert\widehat{\Phi}\vert^{2} + \rho_3\vert\widehat{\Theta}\vert^{2} +k_0 \vert\widehat{\mathbf p}\vert^{2} \right](\xi,\,t)\nonumber\\
&&+\frac{{\tau}_0}{2} \xi^4\displaystyle\int_{0}^{+\infty}\,g(s)\vert\widehat{\Lambda}(\xi,\,t)\vert^{2} \,ds,  
\end{eqnarray}  
where $\tau_0$ is defined by \eqref{tau0}. Applying the operator $\partial_{\xi}$ to \eqref{SG12}, we obtain the system
\begin{align}\label{e101part}
\begin{cases}
\widehat{\mathbf v}_{t} - i\xi\widehat{\mathbf u} = i\widehat{u}, \\
\rho_1 \widehat{\mathbf u}_{t} - ik_1 \xi \,\widehat{\mathbf v} +k_0 \widehat{\mathbf p}= ik_1 \widehat{v}, \\
\widehat{\mathbf z}_{t} - i\xi\widehat{\mathbf y} = i\widehat{y}, \\
\rho_2\widehat{\mathbf y}_{t} - ik_2\xi\,\widehat{\mathbf z} + k_0 \widehat{\mathbf p} = ik_2 \widehat{z}, \\
\widehat{\Phi}_{t} +\xi^2\widehat{\Theta} = -2\xi \widehat{\theta},  \\
\rho_3\widehat{\Theta}_{t} -(k_3 -\tau_0 g_0 )\xi^2 \widehat{\Phi}- ilk_0\,\xi\widehat{\mathbf p} + {\mathbf G} = G_0 ,\\
\widehat{\mathbf p}_{t} -\widehat{\mathbf u} -\widehat{\mathbf y} -il\xi \widehat{\Theta} =il\widehat{\theta},\\
\widehat{\Lambda}_{t} +\widehat{\Lambda}_s -\widehat{\Theta}=0,
\end{cases}
\end{align}
where 
\begin{equation}\label{G12part}
{\mathbf G} = (1-\tau_0 )\gamma \widehat{\Theta}+\tau_0 \xi^4 \int_0^{+\infty} g(s)\widehat{\Lambda}ds\,\,\hbox{and}\,\, G_0 =2(k_3 -\tau_0 g_0 ) \xi \widehat{\phi}+ilk_0\widehat{p}-4\tau_0 \xi^3 \displaystyle\int_{0}^{+\infty}\,g(s){\widehat{\eta}} (\xi ,s) \,ds.
\end{equation}
As for the proof of \eqref{Ep21}, multiplying \eqref{e101part}$_1$-\eqref{e101part}$_7$ by $k_1 \bar{\widehat{\mathbf v}}$, 
$\bar{\widehat{\mathbf u}}$, $k_2 \bar{\widehat{\mathbf z}}$, 
$\bar{\widehat{\mathbf y}}$, $(k_3 - \tau_0 g_0 )\bar{\widehat{\Phi}}$, 
$\bar{\widehat{\Theta}}$ and $k_0\bar{\widehat{\mathbf p}}$, respectively, multiplying \eqref{e101part}$_8$ by $\tau_0\xi^4 g(s)\bar{\widehat{\Lambda}}$ and integrating on $\mathbb{R}_+$ with respect to $s$, adding all the obtained equations, taking the real part of the resulting expression and using \eqref{hd1}, we have
\begin{eqnarray}\label{Ep21part}
\frac{d}{dt} \widehat{\mathbf E}(\xi,\,t) &=& -(1-\tau_0 )\gamma\,\vert\widehat{\Theta}\vert^{2} +\frac{\tau_0}{2} \xi^4\displaystyle\int_{0}^{+\infty}\,g^{\prime}(s)\vert\widehat{\Lambda}\vert^{2} \,ds\\
&& +Re\left[ik_1 \widehat{u}\bar{\widehat{\mathbf v}}+ik_1 \widehat{v}\bar{\widehat{\mathbf u}}+i k_2\widehat{y}\bar{\widehat{\mathbf z}}+ik_2 \widehat{z}\bar{\widehat{\mathbf y}}-2(k_3 -\tau_0 g_0 )\xi\widehat{\theta}\bar{\widehat{\Phi}}+G_0 \bar{\widehat{\Theta}}+ilk_0\widehat{\theta}\bar{\widehat{\mathbf p}}\right].\nonumber
\end{eqnarray} 
This identity shows that $\widehat{\mathbf E}$ is not necessarily nonincreasing with respect to $t$. Because
\begin{equation*}
\vert\widehat{\mathbf U}(\xi,\,t)\vert^2 = \left[\vert\widehat{\mathbf v}\vert^2 + \vert\widehat{\mathbf u}\vert^{2} + \vert\widehat{\mathbf z}\vert^{2} +\vert\widehat{\mathbf y}\vert^{2} + \vert\widehat{\Phi}\vert^{2} + \vert\widehat{\Theta}\vert^{2} +\vert\widehat{\mathbf p}\vert^{2} \right](\xi,\,t) +\tau_0 \xi^4\displaystyle\int_{0}^{+\infty}\,g(s)\vert\widehat{\Lambda}(\xi,\,t)\vert^{2} \,ds,
\end{equation*}
then we see that, according to the right inequality in \eqref{g2}, 
\begin{equation}\label{equivpart}
\vert\widehat{\mathbf U}\vert^2 \sim \widehat{\mathbf E}. 
\end{equation} 
\vskip0,1truecm
\begin{theorem}\label{lemma32part}
Assume that \eqref{taukrho2} is satisfied. Let $\widehat{U}$ be the solution of \eqref{g6}. Then there exist ${\mathbf c},\,\widetilde{{\mathbf c}}>0$ such that, for any $t\in \mathbb{R}_+$ and for any $\xi\in \mathbb{R}^*$, 
\begin{equation}\label{11part}
\vert\widehat{\mathbf U}(\xi,\,t)\vert \leq
\widetilde{{\mathbf c}}e^{-\,{\mathbf c}\,{f} (\xi)\,t}\left[\vert\widehat{\mathbf U}_0 (\xi)\vert+\left(\xi^{7-2\tau_0} +\xi^{-(4+2\tau_0)}\right)\,\vert\widehat{U}_0 (\xi)\vert\right],
\end{equation}
where $f$ is defined in \eqref{22}.
\end{theorem}
\vskip0,1truecm
\begin{proof}
We observe that the left hand sides of \eqref{e101part} are identical to the ones of \eqref{SG12} if we replace $\widehat{v}$, 
$\widehat{u},\,\widehat{z},\,\widehat{y}$, $\widehat{\phi}$, 
$\widehat{\theta}$, $\widehat{p}$, $\widehat{\eta}$ and $G$ by 
$\widehat{\mathbf v}$, $\widehat{\mathbf u},\,\widehat{\mathbf z},\,\widehat{\mathbf y}$, $\widehat{\Phi}$, $\widehat{\Theta}$, 
$\widehat{\mathbf p}$, $\widehat{\Lambda}$ and ${\mathbf G}$, respectively. So, first, we use the same arguments to get similar differential identities to \eqref{eq31}-\eqref{eq6342}, and second, we treat the additional terms generated by the ones in the right hand sindes of \eqref{e101part}.
\vskip0,1truecm
Multiplying \eqref{e101part}$_{4}$ and \eqref{e101part}$_{3}$ by $i\,\xi\,\overline{\widehat{\mathbf z}}$ and $-i\rho_2\,\xi\,\overline{\widehat{\mathbf y}}$, respectively, adding the resulting equations, taking the real part and using \eqref{hd1}, we obtain
\begin{equation}\label{eq31part}
\frac{d}{dt}Re\left(i\rho_2\,\xi\,\widehat{\mathbf y}\, \overline{\widehat{\mathbf z}}\right) = \xi^{2}\left(\rho_2\vert\widehat{\mathbf y}\vert^{2} -k_2 \vert\widehat{\mathbf z}\vert^{2}\right) + k_0\,Re\left(i\,\xi\,\widehat{\mathbf z}\,\overline{\widehat{\mathbf p}}\right) +Re\left(\rho_2\,\xi\widehat{y}\,\overline{\widehat{\mathbf y}}-k_2\,\xi\widehat{z}\,\overline{\widehat{\mathbf z}} \right).
\end{equation}
Similarily, multiplying \eqref{e101part}$_{2}$ and \eqref{e101part}$_{1}$ by $i\,\xi\,\overline{\widehat{\mathbf v}}$ and $-i\rho_1\,\xi\,\overline{\widehat{\mathbf u}}$, respectively, adding the resulting equations, taking the real part and using \eqref{hd1}, we find
\begin{equation}\label{eq32part}
\frac{d}{dt}Re\left(i\rho_1\,\xi\,\widehat{\mathbf u}\, \overline{\widehat{\mathbf v}}\right) = \xi^{2}\left(\rho_1 \vert\widehat{\mathbf u}\vert^{2} -k_1 \vert\widehat{\mathbf v}\vert^{2}\right) + k_0 Re\left(i\,\xi\,\widehat{\mathbf v}\,\overline{\widehat{\mathbf p}}\right) +Re\left(\rho_1\,\xi\widehat{u}\,\overline{\widehat{\mathbf u}}-k_1\,\xi\widehat{v}\,\overline{\widehat{\mathbf v}} \right).
\end{equation}
Also, multiplying \eqref{e101part}$_{6}$ and \eqref{e101part}$_{5}$ by 
$-\overline{\widehat{\Phi}}$ and $-\rho_3\,\overline{\widehat{\Theta}}$, respectively, adding the resulting equations, taking the real part and using \eqref{hd1}, we get
\begin{eqnarray}\label{eq33part}
\frac{d}{dt}Re\left(-\rho_3\,\widehat{\Theta}\, \overline{\widehat{\Phi}}\right) &=& \xi^{2}\left(\rho_3\vert\widehat{\Theta}\vert^{2} -(k_3 -\tau_0 g_0 ) \vert\widehat{\Phi}\vert^{2}  \right) +lk_0\, Re\left(i\,\xi\,\widehat{\Phi}\,\overline{\widehat{\mathbf p}}\right) \\
&& +\,Re\left( \overline{\widehat{\Phi}}\,{\mathbf G}\right)+Re\left(2\rho_3\xi\widehat{\theta}\,\overline{\widehat{\Theta}}-G_0\,\overline{\widehat{\Phi}} \right).\nonumber
\end{eqnarray}
Multiplying \eqref{e101part}$_{4}$ and \eqref{e101part}$_{7}$ by $\xi^2\,\overline{\widehat{\mathbf p}}$ and $\rho_2\xi^2\,\overline{\widehat{\mathbf y}}$, respectively, adding the resulting equations, taking the real part and using \eqref{hd1}, we obtain
\begin{eqnarray}\label{Sandeq35part}
\frac{d}{dt}Re\left(\rho_2\xi^2\,\widehat{\mathbf p}\, \overline{\widehat{\mathbf y}}\right) &=& \xi^{2}\left(\rho_2\,\vert\widehat{\mathbf y}\vert^{2} - k_0\vert\widehat{\mathbf p}\vert^{2}\right) + \rho_2 \xi^{2}\,Re\left(\widehat{\mathbf y}\,\overline{\widehat{\mathbf u}}\right)+ l\rho_2 \xi^{2}\,Re\left(i\xi\,\widehat{\Theta}\,\overline{\widehat{\mathbf y}}\right)  \\
&& + k_2\xi^{2}\,Re\left(i\,\xi\,\widehat{\mathbf z}\,\overline{\widehat{\mathbf p}}\right)+Re\left(ik_2 \xi^2\widehat{z}\,\overline{\widehat{\mathbf p}}+il\rho_2 \,\xi^2\widehat{\theta}\,\overline{\widehat{\mathbf y}} \right).\nonumber
\end{eqnarray}
After, multiplying \eqref{e101part}$_{3}$ and \eqref{e101part}$_{6}$ by $\rho_3\,\overline{\widehat{\Theta}}$ and $\overline{\widehat{\mathbf z}}$, respectively, adding the resulting equations, taking the real part and using \eqref{hd1}, we entail
\begin{eqnarray}\label{eq36part}
\frac{d}{dt}Re\left(\rho_3\widehat{\mathbf z}\, \overline{\widehat{\Theta}}\right) &=& -\rho_3\, Re\left(i\xi\widehat{\Theta}\,\overline{\widehat{\mathbf y}}\right)+(k_3 -\tau_0 g_0)\,\xi^{2}\, Re\left(\widehat{\Phi}\,\overline{\widehat{\mathbf z}}\right)-lk_0\,Re\left(i\,\xi\,\widehat{\mathbf z}\,\overline{\widehat{\mathbf p}}\right)\\
&& -\,Re\left(\overline{\widehat{\mathbf z}}\, {\mathbf G}\right)+Re\left(i\rho_3 \widehat{y}\,\overline{\widehat{\Theta}}+ G_0\,\overline{\widehat{\mathbf z}} \right).\nonumber
\end{eqnarray}
Multiplying \eqref{e101part}$_{5}$ and \eqref{e101part}$_{4}$ by $i\rho_2 \xi\,\overline{\widehat{\mathbf y}}$ and $-i \xi\overline{\widehat{\Phi}}$, respectively, adding the resulting equations, taking the real part and using \eqref{hd1}, it follows that
\begin{eqnarray}\label{eq37part}
\frac{d}{dt}Re\left(i\rho_2 \xi\,\widehat{\Phi}\, \overline{\widehat{\mathbf y}}\right) &=& -\rho_2\xi^{2}\, Re\left(i \xi\widehat{\Theta}\,\overline{\widehat{\mathbf y}}\right)+k_2 \xi^2\, Re\left(\widehat{\Phi}\,\overline{\widehat{\mathbf z}}\right)-k_0\,Re\left(i \xi\widehat{\Phi}\,\overline{\widehat{\mathbf p}}\right)\\
&& +Re\left(k_2\,\xi\widehat{z}\,\overline{\widehat{\Phi}} -2i\rho_2\,\xi^2\widehat{\theta}\,\overline{\widehat{\mathbf y}}\right).\nonumber
\end{eqnarray}
Next, multiplying \eqref{e101part}$_{2}$ and \eqref{e101part}$_{3}$ by $i \xi\overline{\widehat{\mathbf z}}$ and $-i\rho_1 \xi\,\overline{\widehat{\mathbf u}}$, respectively, adding the resulting equations, taking the real part and using \eqref{hd1}, it appears that
\begin{equation}\label{eq310part}
\frac{d}{dt}Re\left(i\rho_1 \xi\,\widehat{\mathbf u}\, \overline{\widehat{\mathbf z}}\right) = -k_1 \xi^2\, Re\left(\widehat{\mathbf v}\,\overline{\widehat{\mathbf z}}\right) +\rho_1  \xi^2\,Re\left(\widehat{\mathbf y} \,\overline{\widehat{\mathbf u}}\right)+k_0 \,Re\left(i \xi\widehat{\mathbf z}\, \overline{\widehat{\mathbf p}}\right)+Re\left(\rho_1\,\xi\widehat{y}\,\overline{\widehat{\mathbf u}} -k_1\,\xi\widehat{v}\,\overline{\widehat{\mathbf z}}\right).
\end{equation}
Multiplying \eqref{e101part}$_{2}$ and \eqref{e101part}$_{5}$ by $i \xi\overline{\widehat{\Phi}}$ and $-i\rho_1 \xi\,\overline{\widehat{\mathbf u}}$, respectively, adding the resulting equations, taking the real part and using \eqref{hd1}, we see that
\begin{eqnarray}\label{eq312part}
\frac{d}{dt}Re\left(i\rho_1 \xi\,\widehat{\mathbf u}\, \overline{\widehat{\Phi}}\right) &=& -k_1 \xi^2\,Re\left(\widehat{\Phi}\,\overline{\widehat{\mathbf v}}\right) +\rho_1 \,\xi^2\,Re\left(i \xi\widehat{\Theta} \,\overline{\widehat{\mathbf u}}\right)+k_0\,Re\left(i \xi\widehat{\Phi}\, \overline{\widehat{\mathbf p}}\right)\\
&& +Re\left(2i\rho_1\,\xi^2\widehat{\theta}\,\overline{\widehat{\mathbf u}} -k_1\,\xi\widehat{v}\,\overline{\widehat{\Phi}}\right).\nonumber
\end{eqnarray}
Also, multiplying \eqref{e101part}$_{1}$ and \eqref{e101part}$_{4}$ by $i\rho_2 \xi\,\overline{\widehat{\mathbf y}}$ and $-i \xi\overline{\widehat{\mathbf v}}$, respectively, adding the resulting equations, taking the real part and using \eqref{hd1}, we find
\begin{equation}\label{eq36part}
\frac{d}{dt}Re\left(i\rho_2 \xi\widehat{\mathbf v}\, \overline{\widehat{\mathbf y}}\right) =- \rho_2 \xi^2\, Re\left(\widehat{\mathbf y}\,\overline{\widehat{\mathbf u}}\right)+k_2 \xi^2\, Re\left(\widehat{\mathbf v}\,\overline{\widehat{\mathbf z}}\right)-k_0\,Re\left(i \xi\widehat{\mathbf v}\,\overline{\widehat{\mathbf p}}\right)+Re\left(k_2\xi\widehat{z}\,\overline{\widehat{\mathbf v}} -\rho_2\,\xi\widehat{u}\,\overline{\widehat{\mathbf y}}\right).
\end{equation}
Multiplying \eqref{e101part}$_{1}$ and \eqref{e101part}$_{6}$ by $\rho_3\,\overline{\widehat{\Theta}}$ and $\overline{\widehat{\mathbf v}}$, respectively, adding the resulting equations, taking the real part and using \eqref{hd1}, we get
\begin{eqnarray}\label{eq37part}
\frac{d}{dt}Re\left(\rho_3\,\widehat{\mathbf v}\, \overline{\widehat{\Theta}}\right) &=& -\rho_3\, Re\left(i\xi\widehat{\Theta}\,\overline{\widehat{\mathbf u}}\right)+(k_3 -\tau_0 g_0 )\, \xi^2\,Re\left(\widehat{\Phi}\,\overline{\widehat{\mathbf v}}\right)-lk_0\,Re\left(i\xi\widehat{\mathbf v}\,\overline{\widehat{\mathbf p}}\right)\\
&& -\,Re\left(\overline{\widehat{\mathbf v}}\, {\mathbf G}\right)+Re\left(i\rho_3\widehat{u}\,\overline{\widehat{\Theta}} +G_0\,\overline{\widehat{\mathbf v}}\right).\nonumber
\end{eqnarray}
Now, multiplying \eqref{e101part}$_{6}$ by $-\xi^2 \displaystyle\int_{0}^{+\infty}\,g(s)\overline{\widehat{\Lambda}} \,ds$, multiplying \eqref{e101part}$_{8}$ by $-\rho_3 \xi^2 g(s)\overline{\widehat{\Theta}}$ and integrating over $\mathbb{R}_+$ with respect to $s$, adding the resulting equations, taking the real part and using \eqref{hd1}, we infer that 
\begin{eqnarray}\label{eq6341part}
\frac{d}{dt}Re\left(-\rho_3 \xi^2 {\widehat{\Theta}}\, \displaystyle\int_{0}^{+\infty}\,g(s)\overline{\widehat{\Lambda}} \,ds\right) &=& -\rho_3 g_0 \xi^2 \vert \widehat{\Theta}\vert^2  -(k_3 -\tau_0 g_0)\xi^4 Re\left(\widehat{\Phi}\displaystyle\int_{0}^{+\infty}\,g(s)\overline{\widehat{\Lambda}} \,ds\right) \nonumber \\ 
&&-lk_0 \xi^2 Re\left(i\xi\widehat{\mathbf p}\displaystyle\int_{0}^{+\infty}\,g(s)\overline{\widehat{\Lambda}} \,ds\right) -\rho_3 \xi^2  Re\left(\overline{\widehat{\Theta}}\displaystyle\int_{0}^{+\infty}\,g^{\prime}(s){\widehat{\Lambda}} \,ds\right)\\
&&+\xi^2 Re \left({\mathbf G}\displaystyle\int_{0}^{+\infty}\,g(s)\overline{{\widehat{\Lambda}}} \,ds\right)-\xi^2 Re \left(G_0\displaystyle\int_{0}^{+\infty}\,g(s)\overline{{\widehat{\Lambda}}} \,ds\right). \nonumber 
\end{eqnarray}
Similarily, multiplying \eqref{e101part}$_{4}$ by 
$i\xi \displaystyle\int_{0}^{+\infty}\,g(s)\overline{\widehat{\Lambda}} \,ds$, multiplying \eqref{e101part}$_{8}$ by $-i\rho_2 \xi g(s)\,\overline{\widehat{\mathbf y}}$ and integrating over $\mathbb{R}_+$ with respect to $s$, adding the resulting equations, taking the real part and using \eqref{hd1}, we arrive at
\begin{eqnarray}\label{Sandeq634part}
\frac{d}{dt}Re\left(i\rho_2 \xi {\widehat{\mathbf y}}\, \displaystyle\int_{0}^{+\infty}\,g(s)\overline{\widehat{\Lambda}} \,ds\right) &=& -\rho_2 g_0 Re\left(i\xi{\widehat{\Theta}} \overline{\widehat{\mathbf y}}\right) -k_2 \xi^2 Re\left(\widehat{\mathbf z}\displaystyle\int_{0}^{+\infty}\,g(s)\overline{\widehat{\Lambda}} \,ds\right) \nonumber \\ 
&&-k_0 Re\left(i\xi\widehat{\mathbf p}\displaystyle\int_{0}^{+\infty}\,g(s)\overline{\widehat{\Lambda}} \,ds\right) -\rho_2 Re\left(i\xi\overline{\widehat{\mathbf y}}\displaystyle\int_{0}^{+\infty}\,g^{\prime}(s){\widehat{\Lambda}} \,ds\right)\\
&& -k_2\,Re \left(\xi\,\widehat{z}\displaystyle\int_{0}^{+\infty}\,g(s)\overline{{\widehat{\Lambda}}} \,ds\right). \nonumber 
\end{eqnarray}
Finally, multiplying \eqref{e101part}$_{2}$ by $i \xi\displaystyle\int_{0}^{+\infty}\,g(s)\overline{\widehat{\Lambda}} \,ds$, multiplying \eqref{e101part}$_{8}$ by $-\rho_1 \xi g(s)\overline{\widehat{\mathbf u}}$ and integrating over $\mathbb{R}_+$ with respect to $s$, adding the resulting equations, taking the real part and using \eqref{hd1}, we see that
\begin{eqnarray}\label{eq6342part}
\frac{d}{dt}Re\left(i\rho_1 \xi {\widehat{\mathbf u}}\, \displaystyle\int_{0}^{+\infty}\,g(s)\overline{{\widehat{\Lambda}}}\,ds\right) &=& -k_1 \xi^2 Re\left(\widehat{\mathbf v}\displaystyle\int_{0}^{+\infty}\,g(s)\overline{\widehat{\Lambda}} \,ds\right)-k_0 Re\left(i \xi\widehat{\mathbf p}\displaystyle\int_{0}^{+\infty}\,g(s)\overline{\widehat{\Lambda}} \,ds\right) \nonumber \\ 
&& -\rho_1 Re\left(i \xi\overline{\widehat{\mathbf u}}\displaystyle\int_{0}^{+\infty}\,g^{\prime}(s){\widehat{\Lambda}} \,ds\right)-\rho_1 g_0 Re\left(i \xi{\widehat{\Theta}}\overline{\widehat{\mathbf u}}\right) \\
&& -k_1\,Re \left(\xi\,\widehat{v}\displaystyle\int_{0}^{+\infty}\,g(s)\overline{{\widehat{\Lambda}}} \,ds\right). \nonumber 
\end{eqnarray}
\vskip0,1truecm
Let ${\mathbf F}_0 ,\cdots,{\mathbf F}_5$ as ${F}_0 ,\cdots,{F}_5$, respectively, with $\widehat{\mathbf v}$, $\widehat{\mathbf u},\,\widehat{\mathbf z},\,\widehat{\mathbf y}$, $\widehat{\Phi}$, 
$\widehat{\Theta}$, $\widehat{\mathbf p}$, $\widehat{\Lambda}$ and 
${\mathbf G}$ instead of $\widehat{v}$, 
$\widehat{u},\,\widehat{z},\,\widehat{y}$, $\widehat{\phi}$, 
$\widehat{\theta}$, $\widehat{p}$, $\widehat{\eta}$ and $G$, respectively. Exploiting \eqref{eq31part}-\eqref{eq6342part}, we find (instead of \eqref{dtF1})
\begin{equation}\label{dtF1part}
\frac{d}{dt}{\mathbf F}_0 (\xi,\,t)={\mathbf F}_5 (\xi,\,t) +{\mathbf F}_4 (\xi,\,t)+{\mathbf F}_3 (\xi,\,t) +R(\xi,t), 
\end{equation} 
where $R$ gathers the mixed products of the components of $U$ and 
${\mathbf U}$, and it is given by
\begin{eqnarray}\label{Rpart} 
R(\xi,\,t) &=& Re\left[-k_2 \lambda_1\xi
\widehat{z}\bar{\widehat{\mathbf z}} +\rho_2 \lambda_1\xi
\widehat{y}\bar{\widehat{\mathbf y}}-k_1 \lambda_2\xi
\widehat{v}\bar{\widehat{\mathbf v}}+\rho_1 \lambda_2\xi
\widehat{u}\bar{\widehat{\mathbf u}}-\lambda_3
G_0 \bar{\widehat{\Phi}} +2\rho_3 \lambda_3\xi
\widehat{\theta}\bar{\widehat{\Theta}}\right]\nonumber \\
&& +Re\left[ ik_2 \lambda_4\xi^2
\widehat{z}\bar{\widehat{\mathbf p}}+i\rho_2 l\lambda_4\xi^2
\widehat{\theta}\bar{\widehat{\mathbf y}}+i\rho_3 \lambda_5
\widehat{y}\bar{\widehat{\Theta}}+\lambda_5 {G_0}\bar{\widehat{\mathbf z}} -2i\rho_2 \lambda_6\xi^2\widehat{\theta}\bar{\widehat{\mathbf y}}+k_2 \lambda_6\xi\widehat{z}\bar{\widehat{\Phi}}\right]\\
&& +Re\left[-k_1 \lambda_7\xi\widehat{v}\bar{\widehat{\mathbf z}} +\rho_1 \lambda_7\xi\widehat{y}\bar{\widehat{\mathbf u}}-k_1 \lambda_8\xi\widehat{v}\bar{\widehat{\Phi}}+2i\rho_1\lambda_8\xi^2\widehat{\theta}\bar{\widehat{\mathbf u}}-\rho_2\lambda_9\xi\widehat{u}\bar{\widehat{\mathbf y}} +k_2\lambda_9\xi\widehat{z}\bar{\widehat{\mathbf v}} \right]\nonumber \\
&& +Re\left[i\rho_3\lambda_{10}\widehat{u}\bar{\widehat{\Theta}}+\lambda_{10}G_0\bar{\widehat{\mathbf v}}+\tau_0 \left(-\lambda_{11}\xi^2 G_0 -k_2\lambda_{12}\xi\widehat{z}-k_1\lambda_{13}\xi\widehat{v}\right)\int_0^{+\infty} g(s)\bar{\widehat{\Lambda}}ds\right].
\end{eqnarray}
Considering the same choices of $\lambda_1 \,\cdots,\,\lambda_{13}$ and using \eqref{hd41} and \eqref{hd42} for 
${\tilde\eta}\in\left\{{\widehat{\Lambda}}, \overline{{\widehat{\Lambda}}}\right\}$, we get (instead of \eqref{xiF04100} and \eqref{xiF0414})  
\begin{eqnarray}\label{xiF0414part}
\frac{d}{dt}\left[\xi^{2+2\tau_0} {\mathbf F}_0 (\xi,\,t)\right] &\leq& -c_1 \xi^{4+2\tau_0} \widehat{\mathbf E}(\xi,\,t)+\xi^{2+2\tau_0}\,R(\xi,t)\\
&& +c_2\xi^{2\tau_0}\left(\xi^{10} +1\right)\left[(1-\tau_0 )\vert \widehat{\Theta}\vert^2 -\tau_0 \xi^2\displaystyle\int_{0}^{+\infty}\,g^{\prime}(s)\vert{\widehat{\Lambda}} (\xi ,s)\vert^2 \,ds\right].\nonumber 
\end{eqnarray}  
Using Young's inequality and exploiting \eqref{m123}, \eqref{vertlambdajleq} and \eqref{vertlambdajleq+}, we obtain, for any 
$\epsilon >0$,
\begin{eqnarray*}
\xi^{2+2\tau_0}\vert R(\xi,\,t)\vert &\leq& \epsilon \xi^{4+2\tau_0} \left(\vert{\widehat{\mathbf z}}\vert^2 +\vert{\widehat{\mathbf y}}\vert^2 +\vert{\widehat{\mathbf v}}\vert^2 +\vert{\widehat{\mathbf u}}\vert^2 +\vert{\widehat{\Phi}}\vert^2 +\vert{\widehat{\Theta}}\vert^2 +\vert{\widehat{\mathbf p}}\vert^2 \right)\\
&& +\tau_0 \epsilon\xi^{6+4\tau_0}\displaystyle\int_{0}^{+\infty}\,g(s)\vert{\widehat{\Lambda}} (\xi ,s)\vert^2 \,ds +C_{\epsilon} \xi^{2\tau_0}\left(\xi^{8} +1\right)\vert{\widehat{U}}\vert^2 ,
\end{eqnarray*}
thus, using \eqref{equivpart}, we obtain in both cases $\tau_0 =0$ and $\tau_0 =1$ 
\begin{equation}\label{R1}
\xi^{2+2\tau_0}\vert R(\xi,\,t)\vert\leq \epsilon C_1 \xi^{4+2\tau_0}{\widehat{\mathbf E}}(\xi,\,t)+C_{\epsilon}\xi^{2\tau_0}\left(\xi^{8}+1\right)\vert\widehat{U}(\xi,\,t)\vert^2 ,
\end{equation}
therefore, chosing $\epsilon =\frac{c_1}{2C_1}$, we deduce from \eqref{xiF0414part} and \eqref{R1} that
\begin{eqnarray}\label{R1xiF0}
\frac{d}{dt} \left[\xi^{2+2\tau_0} {\mathbf F}_0 (\xi,\,t)\right] &\leq& -\frac{c_1}{2}\,\xi^{4+2\tau_0}\,{\widehat{\mathbf E}} (\xi,\,t) +C\xi^{2\tau_0}\left(\xi^{8}+1\right)\vert\widehat{U}(\xi,\,t)\vert^2 \\
&& +c_2 \xi^{2\tau_0}\left(\xi^{10} +1\right)\left[(1-\tau_0 )\vert \widehat{\Theta}\vert^2 -\tau_0 \xi^2\displaystyle\int_{0}^{+\infty}\,g^{\prime}(s)\vert{\widehat{\Lambda}} (\xi ,s)\vert^2 \,ds\right].\nonumber
\end{eqnarray} 
Now, let $\lambda >0$ and ${\mathbf F}$ defined as $F$ in \eqref{xiL1100} and \eqref{xiL1144}; that is 
\begin{equation*}
{\mathbf F} (\xi,\,t)=\lambda\,{\widehat{\mathbf E}} (\xi,\,t) +\frac{\xi^{2+2\tau_0}}{\xi^{10}+1}{\mathbf F}_0 (\xi,\,t).
\end{equation*}
By combining \eqref{Ep21part} and \eqref{R1xiF0}, we arrive at 
\begin{eqnarray}\label{R1xiF0EE}
\frac{d}{dt} {\mathbf F} (\xi,\,t) &\leq& -\frac{c_1 \xi^{4+2\tau_0}}{2\left(\xi^{10}+1\right)}\,{\widehat{\mathbf E}} (\xi,\,t)+C\frac{\xi^{2\tau_0}\left(\xi^{8}+1\right)}{\xi^{10}+1}\vert\widehat{U}(\xi,\,t)\vert^2 +R_1 (\xi,\,t) \nonumber \\
&& +(1-\tau_0)\left(c_2 \xi^{2\tau_0}-\gamma\lambda\right)\vert\widehat{\Theta}\vert^2 +\tau_0\left(\frac{1}{2}\lambda\xi^4 -c_2\xi^{2+2\tau_0}\right)\displaystyle\int_{0}^{+\infty}\,g^{\prime}(s)\vert{\widehat{\Lambda}} (\xi ,s)\vert^2 \,ds, 
\end{eqnarray}
where
\begin{equation*}
R_1 (\xi,\,t)=\lambda\,Re\left[ik_1 \widehat{u}\bar{\widehat{\mathbf v}}+ik_1 \widehat{v}\bar{\widehat{\mathbf u}}+i k_2\widehat{y}\bar{\widehat{\mathbf z}}+ik_2 \widehat{z}\bar{\widehat{\mathbf y}}-2(k_3 -\tau_0 g_0 )\xi\widehat{\theta}\bar{\widehat{\Phi}}+G_0 \bar{\widehat{\Theta}}+ilk_0\widehat{\theta}\bar{\widehat{\mathbf p}}\right].
\end{equation*}
On the other hand, as \eqref{xiL3} and \eqref{xiL3000}, we have
\begin{equation}\label{xiL3E} 
\vert{\mathbf F}(\xi,\,t) - \lambda\,\widehat{\mathbf E}(\xi,\,t)\vert\leq 2c_{3}\,\widehat{\mathbf E}(\xi,\,t).
\end{equation}
Hence, for $\lambda>\left\{\frac{c_2}{\gamma},2c_3\right\}$ in case $\tau_0 =0$, and $\lambda>\left\{2c_2,2c_3\right\}$ in case $\tau_0 =1$, we deduce from \eqref{R1xiF0EE} and \eqref{xiL3E} that  
\begin{equation}\label{R1xiF0E}
\frac{d}{dt} {\mathbf F} (\xi,\,t) \leq -\frac{c_1 \xi^{4+2\tau_0}}{2\left(\xi^{10}+1\right)}\,{\widehat{\mathbf E}} (\xi,\,t)+C\frac{\xi^{2\tau_0}\left(\xi^{8}+1\right)}{\xi^{10}+1}\vert\widehat{U}(\xi,\,t)\vert^2 +R_1 (\xi,\,t)
\end{equation}
and
\begin{equation}\label{equivE}
{\mathbf F} \sim {\widehat{\mathbf E}}.
\end{equation}
Applying Young's inequality and using the definitin of $\widehat{U}$, 
${\widehat{\mathbf E}}$ and $G_0$, we see that, for any $\varepsilon >0$,
\begin{equation}\label{R1xiF0EEE}
\xi^{4+2\tau_0}R_1 (\xi,\,t)\leq \frac{\varepsilon\,\xi^{2(4+2\tau_0)}}{\xi^{10}+1}{\widehat{\mathbf E}} (\xi,\,t) +C_{\varepsilon}\left(\xi^2 +1\right)\left(\xi^{10} +1\right) \vert\widehat{U}(\xi,\,t)\vert^2 ,
\end{equation}  
then, by multiplying \eqref{R1xiF0E} by $\xi^{4+2\tau_0}$, combining with \eqref{R1xiF0EEE}, choosing $\varepsilon=\frac{c_1}{4}$ and using \eqref{m123} and \eqref{equivE}, we obtain, for some $c_0 =\frac{c_1}{8}$, 
\begin{equation}\label{R1xiF0Epart}
\frac{d}{dt}\left[\xi^{4+2\tau_0}{\mathbf F} (\xi,\,t)\right] \leq -\frac{2c_0 \xi^{2(4+2\tau_0)}}{\xi^{10}+1}\,{\mathbf F} (\xi,\,t)+C\frac{\xi^{22}+1}{\xi^{10}+1}\vert\widehat{U}(\xi,\,t)\vert^2 , 
\end{equation}
therefore, by multiplying \eqref{R1xiF0Epart} by $e^{2c_0 f(\xi)t}$ ($f$ is defined in \eqref{22}) and using \eqref{11}, we find   
\begin{equation*}
\frac{d}{dt} \left[\xi^{4+2\tau_0}e^{2c_0 f(\xi)t}{\mathbf F} (\xi,\,t)\right] \leq C\frac{\xi^{22} +1}{\xi^{10}+1}e^{2(c_0 -c) f(\xi)t}\vert\widehat{U}_0(\xi)\vert^2 ,
\end{equation*}
then, by integrating the above inequality with respect to $t$, we find
\begin{equation}\label{R2}
{\mathbf F}(\xi,\,t)\leq e^{-2c_0 f(\xi)t}{\mathbf F}(\xi,\,0)+C
\frac{\xi^{22} +1}{\xi^{2(4+2\tau_0 )}}
e^{-2c f(\xi)t}\vert\widehat{U}_0(\xi)\vert^2 ,
\end{equation}
thus, according to \eqref{equivpart} and \eqref{equivE} and using the inequality \eqref{sqrta1a2}, the above inequality \eqref{R2} implies \eqref{11part}. The proof of Theorem \ref{lemma32part} is now achieved.
\end{proof}

\section{Estimation of $\left\Vert\partial_{x}^{j}U\right\Vert_{1}$}

In this section, we show the decay estimate of 
$\Vert\partial_{x}^{j}U\Vert_{1}$ by exploiting \eqref{sandwi332}, \eqref{sandwi334} and \eqref{11part}. 
\vskip0,1truecm
\begin{theorem}\label{theorem41L1}
Assume that \eqref{taukrho2} is satisfied. Let $N\,\in \{5+2\tau_0 ,6+2\tau_0,\cdots\}$ and $U$ be a solution of \eqref{e11} corresponding to an initial data $U_{0}$ satisfying 
\begin{equation}\label{regularityL2}
U_{0}\in H^{N+7-2\tau_0}(\mathbb{R})\cap L^{1}(\mathbb{R})\quad\hbox{and}\quad {\tilde U}_{0}\in H^{N}(\mathbb{R})\cap L^{1}(\mathbb{R}),
\end{equation}
where ${\tilde U}_{0} (x)=x {U}_{0} (x)$. Then, for any $\ell\in \{1,\,2,\,\cdots,\,N-4-2\tau_0\}$ and $j\in\{4+2\tau_0,\,5+2\tau_0,\,\cdots,\,N-\ell\}$, there exist ${\mathbf c}_0 ,\,{\tilde{\mathbf c}}_0 >0$ such that, for any 
$t\in \mathbb{R}_+$, 
\begin{eqnarray}\label{33L1}
\Vert\partial_{x}^{j}U\Vert_{{1}} &\leq& {\tilde {\mathbf c}}_0 \,\left[(1 + t)^{-1/8 - j/4}\,\Vert {\tilde U}_{0} \Vert_{{1}} +(1 + t)^{7/8 - j/4}\,\Vert {U}_{0} \Vert_{{1}}\right]\\
&& + {\tilde {\mathbf c}}_0 \,(1 + t)^{-\ell/6}\,\left[\Vert\partial_{x}^{j+\ell}{\tilde U}_{0} \Vert_{{2}}+\Vert\partial_{x}^{j+\ell+7}U_{0} \Vert_{{2}}+\Vert\partial_{x}^{j+\ell}U_{0} \Vert_{{2}} +\Vert\partial_{x}^{j+\ell-1}U_{0} \Vert_{{2}} \right].\nonumber
\end{eqnarray}
if $\tau_0=0$, and 
\begin{eqnarray}\label{3300L1}
\Vert\partial_{x}^{j}U\Vert_{{1}} &\leq& {\tilde {\mathbf c}}_0 \,\left[(1 + t)^{-1/12 - j/6}\,\Vert {\tilde U}_{0} \Vert_{{1}} +(1 + t)^{11/12 - j/6}\,\Vert {U}_{0} \Vert_{{1}}\right]\\
&& +{\tilde {\mathbf c}}_0 (1 + t)^{-\ell/4}\,\left[\Vert\partial_{x}^{j+\ell}{\tilde U}_{0} \Vert_{{2}}+\Vert\partial_{x}^{j+\ell+5}U_{0} \Vert_{{2}}+\Vert\partial_{x}^{j+\ell}U_{0} \Vert_{{2}} +\Vert\partial_{x}^{j+\ell-1}U_{0} \Vert_{{2}} \right]\nonumber
\end{eqnarray}
if $\tau_0=1$.
\end{theorem}
\vskip0,1truecm
\begin{proof}
First, using \eqref{Carlson} and applying Young's inequality, it follows that 
\begin{equation}\label{L1L2}
\Vert\partial_{x}^{j}U\Vert_{1} \leq {\frac{1}{2}}\Vert\partial_{x}^{j}U\Vert_{2}+{\frac{1}{2}}\Vert x\partial_{x}^{j}U\Vert_{2}.
\end{equation}
Because the term $\Vert\partial_{x}^{j}U\Vert_{2}$ has yet been estimated in 
\eqref{sandwi332} and \eqref{sandwi334}, we have only to estimate the term 
$\Vert x\partial_{x}^{j}U\Vert_{2}$. Using the method in \cite{Ign_Ross_2010} and Plancherel's theorem, we may write 
\begin{eqnarray}\label{L1L20}
\int_{\mathbb{R}}x^{2}|\partial _{x}^{j}U(x,t)|^{2}dx &\leq &C\int_{
\mathbb{R}}\left\vert \partial _{\xi }\left( |\xi |^{j}|\widehat{U}(\xi
,t)|\right) \right\vert ^{2}d\xi \\
&\leq &C\int_{\mathbb{R}}\left( |\xi |^{j-1}|\widehat{U}(\xi ,t)|+|\xi
|^{j}|\partial _{\xi }\widehat{U}(\xi ,t)|\right) ^{2}d\xi \nonumber\\
&\leq &C\Vert \partial _{x}^{j-1}U\Vert _{2}^{2}+C\int_{\mathbb{R}
}\xi^{2j}|\widehat{\mathbf U}(\xi ,t)|^{2}d\xi.\nonumber
\end{eqnarray}
It is clear that $\Vert \partial _{x}^{j-1}U\Vert _{2}$ can be easily estimated by using \eqref{sandwi332} and \eqref{sandwi334} (with $j-1$ instead of $j$). To estimate the last integral in \eqref{L1L20}, we use \eqref{11part} and apply Plancherel's theorem, it appears that
\begin{eqnarray}\label{g37L1}
\int_{\mathbb{R}}\xi^{2\,j}\,\vert\widehat{U}(\xi,\,t)\vert^{2} d\xi 
& \leq & C\int_{\mathbb{R}^*}\xi^{2\,j}\,e^{-2\,{\mathbf c}\,f (\xi)\,t}\,\left[\vert\widehat{\mathbf U}_0 (\xi)\vert^{2} +\left( \xi^{14-4\tau_0}+\xi^{-(8+4\tau_0 )}\right)\vert\widehat{U}_0 (\xi)\vert^{2}\right]\,d\xi \\
& \leq & C\int_{0<\vert\xi\vert\leq 1}\,e^{-2\,{\mathbf c}\,f (\xi)\,t}\,\left[\xi^{2\,j}\vert\widehat{\mathbf U}_0 (\xi)\vert^{2} +\xi^{2(j-4-2\tau_0)}\vert\widehat{U}_0 (\xi)\vert^{2}\right]d\xi \nonumber\\
&& + C\int_{\vert\xi\vert > 1}\,e^{-2\,{\mathbf c}\,f (\xi)\,t}\,\left[\xi^{2j}\vert\widehat{\mathbf U}_0 (\xi)\vert^{2} +\xi^{2\,(j+7-2\tau_0)}\vert\widehat{U}_0 (\xi)\vert^{2}\right]\, d\xi \nonumber\\
& := & J_{1} + J_{2}.\nonumber
\end{eqnarray}
\vskip0,1truecm
{\bf Case $\tau_0 =0$}: using \eqref{Sg36}$_1$ and Lemma 2.3 of \cite{gues2}, we observe that
\begin{eqnarray}\label{g38L1}
J_{1} &\leq& C\,\Vert\widehat{\mathbf U}_{0}\Vert_{{\infty}}^{2}\int_{\vert\xi\vert\leq 1}\xi^{2\,j}\,e^{-\,{\mathbf c}\,t\,\xi^{4}}\, d\xi+
C\,\Vert\widehat{U}_{0}\Vert_{{\infty}}^{2}\int_{\vert\xi\vert\leq 1}\xi^{2\,(j-4)}\,e^{-\,{\mathbf c}\,t\,\xi^{4}}\, d\xi\\
&\leq& C\left(1 + t\right)^{-\,\frac{1}{4}(1 + 2\,j)}\Vert{\tilde U}_{0}\Vert_{{1}}^{2} +C\left(1 + t\right)^{-\,\frac{1}{4}[1 + 2\,(j-4)]}\Vert {U}_{0}\Vert_{{1}}^{2} .\nonumber
\end{eqnarray}
In the high frequency region, using \eqref{Sg36}$_2$, we entail
\begin{eqnarray*}
J_{2} & \leq & C\int_{\vert\xi\vert > 1}\xi^{2\,j}\,e^{-{\mathbf c}\,t\,\xi^{-6}}\,\vert\widehat{\mathbf U}_0 (\xi)\vert^{2}\, d\xi +C\int_{\vert\xi\vert > 1}\xi^{2\,(j+7)}\,e^{-{\mathbf c}\,t\,\xi^{-6}}\,\vert\widehat{U}_0 (\xi)\vert^{2}\, d\xi\\
& \leq & C\ \sup_{\vert\xi\vert >1}\left\{\xi^{-2\,\ell}\, e^{-{\mathbf c}\,t\,\xi^{-6}}\right\} \left[\int_{\mathbb{R}}\xi^{2\,(j + \ell )}\,\vert\widehat{\mathbf U}_0 (\xi)\vert^{2}\ d\xi+\int_{\mathbb{R}}\xi^{2\,(j + \ell +7)}\,\vert\widehat{U}_0 (\xi)\vert^{2}\ d\xi\right]\\
& \leq & C\ \sup_{\vert\xi\vert >1}\left\{\xi^{-2\,\ell}\, e^{-{\mathbf c}\,t\,\xi^{-6}}\right\} \,\left(\Vert\,\partial_{x}^{j + \ell}{\tilde U}_{0}\,\Vert_{{2}}^{2} +\Vert\,\partial_{x}^{j + \ell+7}U_{0}\,\Vert_{{2}}^{2} \right),
\end{eqnarray*}
thus, simple computations (see, for example, Lemma 2.4 of \cite{gues2}) imply that
\begin{equation}\label{g38+L1}
J_{2} \leq C\left(1 + t\right)^{-\frac{\ell}{3}}\,\,\left(\Vert\,\partial_{x}^{j + \ell}{\tilde U}_{0}\,\Vert_{{2}}^{2} +\Vert\,\partial_{x}^{j + \ell+7}U_{0}\,\Vert_{{2}}^{2} \right), 
\end{equation}
and so, by combining \eqref{g37L1}-\eqref{g38+L1}, we get 
\begin{eqnarray}\label{g38+L1LLL}
\int_{\mathbb{R}}\xi^{2\,j}\,\vert\widehat{U}(\xi,\,t)\vert^{2} d\xi &\leq& C\left(1 + t\right)^{-\,\frac{1}{4}(1 + 2\,j)}\Vert{\tilde U}_{0}\Vert_{{1}}^{2} +\left(1 + t\right)^{-\,\frac{1}{4}[1 + 2\,(j-4)]}\Vert {U}_{0}\Vert_{{1}}^{2} \\
&& +C\left(1 + t\right)^{-\frac{\ell}{3}}\,\,\left(\Vert\,\partial_{x}^{j + \ell}{\tilde U}_{0}\,\Vert_{{2}}^{2} +\Vert\,\partial_{x}^{j + \ell+7}U_{0}\,\Vert_{{2}}^{2} \right)\nonumber
\end{eqnarray}
thus, by combining \eqref{sandwi332}, \eqref{L1L2} and \eqref{g38+L1LLL}, and using the inequality \eqref{sqrta1a2}, we get \eqref{33L1}. 
\vskip0,1truecm
{\bf Case $\tau_0 =1$}: using \eqref{sandwi334}, \eqref{Sg360} and \eqref{11part}, and following the same arguments as for the proof of \eqref{33L1}, we obtain \eqref{3300L1}.  
\end{proof}

\section{Concluding discussion}

1. This article is concerned with the stability of two systems of type Rao-Nakra sandwich beam in the whole line $\mathbb{R}$ under the presence of a frictional damping or an infinite memory acting on the Euler-Bernoulli equation. We prove the instability of both systems if the speeds of propagation of the two wave equations are equal. In the reverse situation, we are able, despite the presence of only one controle, to obtain the desired $L^2 (\mathbb{R})$-norm and $L^1 (\mathbb{R})$-norm decay estimates \eqref{sandwi332}, \eqref{sandwi334}, \eqref{33L1} and \eqref{3300L1}. The main ingredient of the proof is the energy method and the Fourier analysis. 
\vskip0,1truecm
2. The decay estimates \eqref{sandwi332}, \eqref{sandwi334}, \eqref{33L1} and \eqref{3300L1} are still satisfied when \eqref{taukrho1} holds if we add 
${\tilde\gamma}\varphi_t$ to \eqref{S1}$_1$ and \eqref{S2}$_1$, or 
${\tilde\gamma}\psi_t$ to \eqref{S1}$_2$ and \eqref{S2}$_2$, where 
${\tilde\gamma}>0$. Indeed, in these cases, we do not need $A_7 =0$ because either $\vert\widehat{u}\vert^{2}$ or $\vert\widehat{y}\vert^{2}$ will be directly controlled via the derivative of the energy functional $\widehat{E}$. Similarily, \eqref{sandwi332}, \eqref{sandwi334}, \eqref{33L1} and \eqref{3300L1} hold even in case \eqref{taukrho1} if the added control is of memory type:   
\begin{equation*}
\displaystyle \int_0^{+\infty} {\tilde g}(s) \, \varphi_{xx} (x,t-s) ds \quad\hbox{or}\quad \displaystyle \int_0^{+\infty} {\tilde g}(s) \, \psi_{xx} (x,t-s) ds   
\end{equation*}
instead of ${\tilde\gamma}\varphi_t$ and ${\tilde\gamma}\psi_t$, respectively, where ${\tilde g}$ is as $g$. 
\vskip0,1truecm
3. Using the interpolation inequalities \eqref{Interpol_1}, \eqref{Interpol_2} and \eqref{Interpol_3}, we see that our $L^2 (\mathbb{R})$-norm and $L^1 (\mathbb{R})$-norm decay estimates lead to similar 
$L^q (\mathbb{R})$-norm ones, for any $q\in [1,+\infty]$. The $L^q$-norm decay estimates, for $1\leq q <2$, are based on the 
$L^1$-norm decay estimate, and so they require initial data $U_0$ having the regularity \eqref{regularityL2} with $N\in\{5,6,\cdots\}$ in case \eqref{S1}, and $N\in\{7,8,\cdots\}$ in case \eqref{S2}. However, the 
$L^q$-norm decay estimates, for $q\in [2,+\infty]$ require the weaker regularity \eqref{U0regularity}, where $N\in\mathbb{N}^*$.
\vskip0,1truecm
4. In a future work, we aspire to treat the case where the control occurs on a one wave equation of systems. On the other hand, we think that similar results can be obtained with a control subject to a thermal effect like Fourier law, Cattaneo law and Gurtin-Pipkin law. These kinds of controls deserve to be treated and we aspire to do it in a future work.   
\vskip0,2truecm
{\bf Acknowledgment}. The author thanks Belkacem Said-Houari for useful and fruitful discussions and exchanges on $L^q (\mathbb{R})$-norm decay estimates for Cauchy PDEs.

\end{document}